\documentclass[12pt]{article}

\usepackage[T1]{fontenc}

\usepackage{amsmath}
\usepackage{amsthm}
\usepackage{amssymb}
\usepackage{latexsym}
\usepackage{nicefrac}

\usepackage{csquotes}
\MakeOuterQuote{"}

\usepackage{xargs}

\usepackage{enumitem}
\usepackage{calc}
\usepackage{comment}

\usepackage{xcolor}

\usepackage{array}
\usepackage{tabularx}

\usepackage{tikz}
\usetikzlibrary{
  matrix,
  calc,
  backgrounds,
  hobby,
  decorations,
  decorations.markings,
  arrows.meta
}

\usepackage[small, pagestyles]{titlesec}
\usepackage[small, it]{caption}

\usepackage[square,comma,numbers,sort&compress]{natbib}

\usepackage[margin=1in,letterpaper,portrait]{geometry}

\theoremstyle{plain}

\newtheorem{theorem}{Theorem}[section]
\newtheorem{proposition}[theorem]{Proposition}

\newtheorem{question}[theorem]{Question}
\newtheorem*{proposition*}{Proposition}

\newtheorem{corollary}[theorem]{Corollary}
\newtheorem{conjecture}[theorem]{Conjecture}
\newtheorem{observation}[theorem]{Observation}

\newtheorem*{theorem*}{Theorem}

\theoremstyle{definition}

\newtheorem{definition}[theorem]{Definition}

\newtheorem*{namedtheoreminner}{\theoremname}

\usepackage{soul}

\newcommand{\graphrefl}{\text{\L}}

\tikzset{every path/.style={
    thick
}}

\tikzset{every node/.style={
    minimum width=5pt,
    inner sep=1pt,
}}

\tikzset{hollow/.style={
    shape=circle,
    fill=white,
    node contents={},
    outer sep={0pt},
    draw
}}

\tikzset{solid/.style={
    shape=circle,
    fill=black,
    node contents={},
    outer sep={0pt},
    draw
}}

\tikzset{
  ->-/.style={
    decoration={
      markings,
      mark=at position 0.6 with {
        \arrow{Latex[length=7pt]}
      }
    },
    postaction={decorate}
  }
}

\tikzset{
  -->-/.style={
    decoration={
      markings,
      mark=at position 0.8 with {
        \arrow{Latex[length=7pt]}
      }
    },
    postaction={decorate}
  }
}

\newcommand{\plotperm}[2][v]{%
  \foreach \j [count=\i] in {#2} {%
    \node[draw] (#1\j) at (\i,\j) [solid] {};
  }%
}

\newcommand{\plotpoints}[2][v]{%
  \foreach \i/\j [count=\k] in {#2} {%
    \node[draw] (#1\k) at (\i,\j) [solid] {};
  }%
}

\newcommandx*{\labelnode}[3][1=45]{%
    \pgfgettransformentries
        \myxscale\myxslant\myyslant\myyscale\myxshift\myyshift
    \node at ([shift=({#1}:{0.1/\myxscale})]#2)
      [anchor={180+#1}] {#3};%
}

\newcommandx*{\drawjoin}[3][1=]{%
    \foreach \nodea in {#2} {%
        \foreach \nodeb in {#3} {%
            \ifx\nodea\nodeb
            \else
                \draw (\nodea)--(\nodeb);
            \fi
        }%
    }%
}

\newcommand{\makeninegon}[3]{
  \foreach \j [count=\i] in {1,2,...,9} {
    \node[draw] (v\i)
      at ({90+(\j-1)*(360/9)}:0.35)
      [solid] {};
  }

  \labelnode[-90]{v5}{#1}
  \labelnode[-90]{v6}{#2}

  \node[] at (0,0) {\normalsize #3};
}

\usepackage[bookmarks]{hyperref}

\hypersetup{
  colorlinks=true,
  linkcolor=black,
  anchorcolor=black,
  citecolor=black,
  urlcolor=black,
  pdfpagemode=UseThumbs,
  pdftitle={Reflections in Graphs},
  pdfsubject={Combinatorics},
  pdfauthor={Sean Mandrick, Bridget Eileen Tenner, and Vincent Vatter},
}

\title{Reflections in graphs}

\author{
  Sean Mandrick%
  \footnote{
    Department of Mathematics,
    University of Florida,
    Gainesville, Florida.
  } \and
  Bridget Eileen Tenner%
  \footnote{
    Department of Mathematical Sciences,
    DePaul University,
    Chicago, Illinois.
    \texttt{bridget@math.depaul.edu}\\
	Research partially supported by NSF Grant DMS-2054436, by a University Research Council Competitive Research Leave from DePaul University, and by an AMS-Simons Research Enhancement Grant for PUI Faculty.
  } 
  \and 
  Vincent Vatter%
  \footnote{
    Department of Mathematics,
    University of Florida,
    Gainesville, Florida.
    \texttt{vatter@ufl.edu}
  }
}

\date{}

\begin{document}

\maketitle


\begin{abstract}
We introduce and study an operation on graphs, called edge reflection, modeled on the effect of a reflection on the inversion graph of a permutation. The edge reflection length of a graph counts the fewest edge reflections needed to reduce it to the edgeless graph. Emptying a graph one edge at a time, or one vertex at a time, gives two upper bounds on this length, and each holds with equality precisely when the graph is a forest. For the edge bound we say more: the cyclomatic number of the graph bounds, above and below, the difference between the graph's size and its edge reflection length. Returning to permutations, we show that the edge reflection length of an inversion graph is at most the reflection length of any permutation realizing it, with equality whenever each connected component of the inversion graph has cyclomatic number at most one. 
\end{abstract}

\section{Introduction}

Permutations, which model so many phenomena throughout mathematics, play a central role in both combinatorics and algebra. In combinatorics they serve as basic models of order, and in algebra they realize the symmetric group $S_n$, which is the finite Coxeter group of type~A. From the algebraic side, the group $S_n$ is generated by its \emph{simple reflections} $s_1, \ldots, s_{n-1}$, where~$s_i$ swaps the values $i$ and $i+1$, and more generally by its \emph{reflections}~$t_{ab}$ with $1 \le a < b \le n$, where~$t_{ab}$ swaps the values~$a$ and~$b$. Two statistics measure how efficiently a permutation can be built from these generators: the \emph{Coxeter length} $\ell(\pi)$ is the minimum number of simple reflections in any expression for~$\pi$, and equals the number of inversions of~$\pi$; the \emph{reflection length} $\ell'(\pi)$ is the minimum number of arbitrary reflections, and equals $n - c(\pi)$, where $c(\pi)$ is the number of cycles of~$\pi$. Equivalently, each length counts the steps needed to \emph{sort}~$\pi$ to the identity, using simple reflections in the first case and arbitrary reflections in the second.

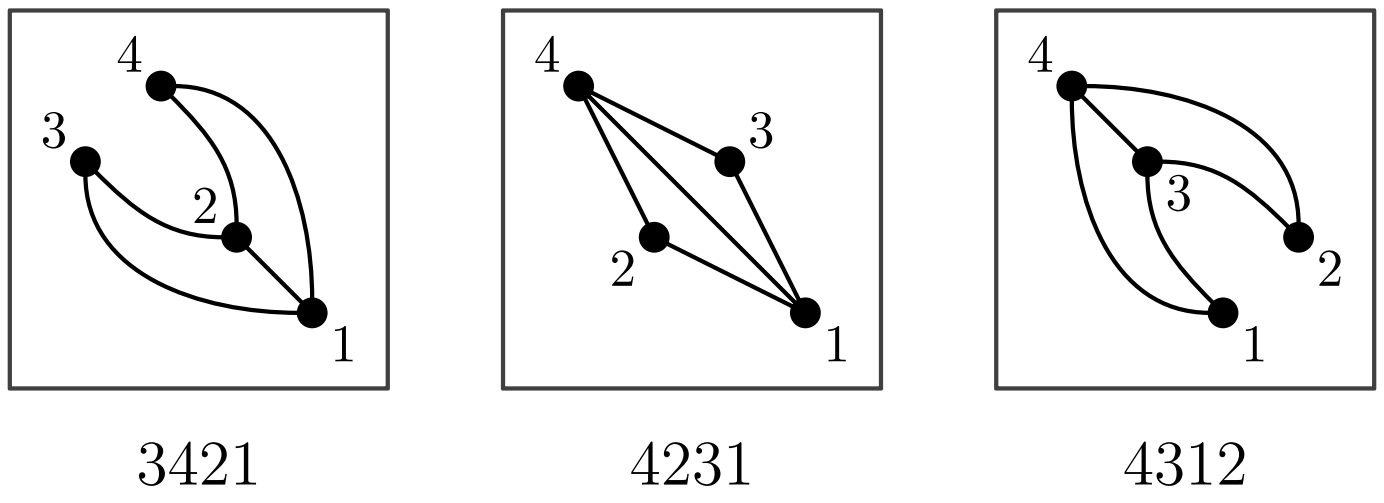
\begin{figure}[htbp]
\begin{footnotesize}
\begin{center}
	\begin{tikzpicture}[scale=0.5]
		\draw[darkgray,thick,rounded corners=0.01,line cap=round] (0,0) rectangle (5,5);
		\plotperm{3,4,2,1}
		\labelnode[135]{v3}{$3$}
		\labelnode[135]{v4}{$4$}
		\labelnode[135]{v2}{$2$}
		\labelnode[-45]{v1}{$1$}
		\draw (v1)--(v2);
		\draw [out=180, in=-90] (v1) to (v3);
		\draw [out=90, in=0] (v1) to (v4);
		\draw [out=180, in=-45] (v2) to (v3);
		\draw [out=90, in=-45] (v2) to (v4);
		\node at (2.5,-1) {\normalsize $3421$};
	\end{tikzpicture}
	\hspace{5mm}
	\begin{tikzpicture}[scale=0.5]
		\draw[darkgray,thick,rounded corners=0.01,line cap=round] (0,0) rectangle (5,5);
		\plotperm{4,2,3,1}
		\labelnode[135]{v4}{$4$}
		\labelnode[-135]{v2}{$2$}
		\labelnode[45]{v3}{$3$}
		\labelnode[-45]{v1}{$1$}
		\draw (v1)--(v2);
		\draw (v1)--(v3);
		\draw (v1)--(v4);
		\draw (v2)--(v4);
		\draw (v3)--(v4);
		\node at (2.5,-1) {\normalsize $4231$};
	\end{tikzpicture}
	\hspace{5mm}
	\begin{tikzpicture}[scale=0.5]
		\draw[darkgray,thick,rounded corners=0.01,line cap=round] (0,0) rectangle (5,5);
		\plotperm{4,3,1,2}
		\labelnode[135]{v4}{$4$}
		\labelnode[-45]{v3}{$3$}
		\labelnode[-45]{v1}{$1$}
		\labelnode[-45]{v2}{$2$}
		\draw [out=135, in=-90](v1) to (v3);
		\draw [out=180, in=-90](v1) to (v4);
		\draw [out=135, in=0] (v2) to (v3);
		\draw [out=90, in=0] (v2) to (v4);
		\draw (v3)--(v4);
		\node at (2.5,-1) {\normalsize $4312$};
	\end{tikzpicture}
\end{center}
\end{footnotesize}
\caption{The permutations $3421$, $4231$, and $4312$, each drawn by plotting the value $\pi(i)$ at the point $(i, \pi(i))$ and joining each inverted pair. All three inversion graphs are isomorphic to $K_4 - e$, yet their reflection lengths are $3$, $1$, and~$3$.}
\label{fig:diamond-realizers}
\end{figure}

A permutation also determines the \emph{inversion graph} $G_\pi$ of $\pi \in S_n$, often called its \emph{permutation graph}. This has vertex set $[n] = \{1, \ldots, n\}$, with an edge joining $a$ and~$b$ whenever these values form an inversion of~$\pi$, meaning that the larger value appears to the left of the smaller in the one-line notation $\pi = \pi(1) \pi(2) \cdots \pi(n)$. Figure~\ref{fig:diamond-realizers} shows three examples, drawing each permutation by placing the value $\pi(i)$ at the point $(i, \pi(i))$ and joining two values whenever they form an inversion. Thus we label the vertices of inversion graphs by \emph{value}, while some other works label them by \emph{position}; the choice is purely cosmetic, but the reader should bear this convention in mind.

Because $\ell(\pi)$ counts inversions, it is a property of the graph alone: $\ell(\pi) = |E(G_\pi)|$. Reflection length, on the other hand, is not a property of the graph alone. The three permutations of Figure~\ref{fig:diamond-realizers} share the inversion graph $K_4 - e$, yet have reflection lengths $3$, $1$, and~$3$, respectively. We return to these three permutations in Section~\ref{sec:permutations}.

Despite these examples, the graph constrains reflection length: every simple reflection is a reflection, so $\ell'(\pi) \le \ell(\pi)$, and the gap $\ell(\pi) - \ell'(\pi)$ counts the steps saved by sorting with arbitrary rather than simple reflections. Tenner and Vatter~\cite{TV} proved that this gap, although not determined by the inversion graph, is bounded by it on both sides:
\begin{equation}
\tag{$\bigstar$}
\label{eqn:tv-gap}
	\beta(G_\pi) \le \ell(\pi) - \ell'(\pi) \le 2\beta(G_\pi),
\end{equation}
where the \emph{cyclomatic number} $\beta(G)$ is the minimum number of edges whose deletion from~$G$ leaves a forest.

This paper lifts the sorting operation itself from inversion graphs to arbitrary graphs. A reflection acts on the inversion graph in a definite way: if $a$ and~$b$ form an inversion of~$\pi$, then applying $t_{ab}$ deletes the edge $ab$ and toggles the adjacencies of $a$ and~$b$ to each value lying horizontally between them in the plot of~$\pi$. An arbitrary graph, on the other hand, has no notion of ``horizontally between,'' so we need another way to choose the vertices that play this intermediate role. Our graph analogue lets that set be chosen freely, requiring only that each of its vertices be adjacent to at least one of the two endpoints (which always holds in the permutation setting). The resulting operation, defined formally in the next section, is an edge reflection, and the edge reflection length $\graphrefl(G)$ counts the fewest edge reflections needed to reduce~$G$ to the edgeless graph.

This freedom makes edge reflections more permissive than the permutation reflections they are modeled on. It is therefore noteworthy how much of the analogy with permutations survives.

We say that a sequence of edge reflections \emph{empties} a graph if it reduces that graph to the edgeless graph. Our first main result concerns the difference $|E(G)| - \graphrefl(G)$, the number of reflections saved over emptying~$G$ one edge at a time. Because $\ell(\pi) = |E(G_\pi)|$, this is the graph analogue of the gap in~\eqref{eqn:tv-gap}, and Theorem~\ref{thm:sandwich} shows that the cyclomatic number bounds it on both sides,
\[
	\beta(G) \le |E(G)| - \graphrefl(G) \le 2\beta(G).
\]
Theorem~\ref{thm:sandwich-improved} improves this lower bound from $\beta(G)$ to $\beta(G) + q(G)$, where $q(G)$ is the number of connected components of~$G$ that contain a cycle.

There are two obvious strategies for emptying a given graph: either empty it edge by edge, or isolate one vertex at a time. Proposition~\ref{prop:upper-bounds} shows that emptying edge by edge gives the bound $\graphrefl(G) \le |E(G)|$, and that isolating vertices gives the bound $\graphrefl(G) \le |V(G)| - \kappa(G)$, where $\kappa(G)$ denotes the number of connected components of~$G$. Forests are the unique extremal graphs for both: they meet the edge bound by Theorem~\ref{thm:sandwich}, since $|E(G)| = \graphrefl(G)$ if and only if $\beta(G) = 0$, and the vertex bound by Theorem~\ref{thm:forest-iff-vertex-bound}.

Returning to permutations, Theorem~\ref{thm:graph-length-at-most-reflection-length} relates edge reflection length and reflection length by $\graphrefl(G_\pi) \le \ell'(\pi)$. Together with~\eqref{eqn:tv-gap} and Theorem~\ref{thm:sandwich}, this gives
\[
	\beta(G_\pi) \le \ell(\pi) - \ell'(\pi) \le \ell(\pi) - \graphrefl(G_\pi) \le 2\beta(G_\pi),
\]
so $\ell'(\pi) - \graphrefl(G_\pi)$ lies between $0$ and $\beta(G_\pi)$. Moreover, as Section~\ref{sec:permutations} shows, this difference is always even, so it must be zero whenever $\beta(G_\pi) \le 1$. Applying this to each connected component of~$G_\pi$, the difference is zero whenever every connected component has cyclomatic number at most one.

\subsection*{Organization}

Section~\ref{sec:edge-reflections} examines the effect of a reflection on an inversion graph, defines edge reflections and edge reflection length, and records the elementary bounds. Section~\ref{sec:sorting} proves that $\graphrefl(G_\pi) \le \ell'(\pi)$. Section~\ref{sec:cyclomatic} proves the cyclomatic bounds and identifies forests as the graphs meeting the edge bound, while Section~\ref{sec:vertex-bound} shows that they are also the graphs meeting the vertex bound and strengthens the lower bound to $\beta(G) + q(G)$. Section~\ref{sec:lower-bounds} gives lower bounds in terms of the vertex cover and edge domination numbers, and computes $\graphrefl$ for several families. Section~\ref{sec:permutations} collects the consequences for permutations, and we close with open questions in Section~\ref{sec:further-questions}.

The results of this paper grew out of the first author's doctoral thesis~\cite{mandrick:on-inversion-gr:}.

\section{Edge reflections}
\label{sec:edge-reflections}

A reflection acts on a permutation by swapping two of its values. We first examine the effect of such a swap on the inversion graph, and then ask how that effect can be described in the graph setting.

\begin{figure}[t]
\begin{footnotesize}
\begin{center}
	\begin{tikzpicture}[scale=0.5]
		\draw [draw = none, fill=gray!20] (2,0.5) rectangle (5,6.5);
		\draw [color = gray, dotted, thick, line cap=round] (2,0.5)--(2,6.5);
		\draw [color = gray, dotted, thick, line cap=round] (5,0.5)--(5,6.5);
		\draw [color = gray, dotted, thick, line cap=round] (0.5,2)--(6.5,2);
		\draw [color = gray, dotted, thick, line cap=round] (0.5,5)--(6.5,5);
		\draw [out=0, in=90] (1,6) to (6,4);
		\draw [out=-45, in=135] (1,6) to (2,5);
		\draw [out=-60, in=165] (1,6) to (4,3);
		\draw [out=-75, in=180] (1,6) to (5,2);
		\draw [out=-90, in=180] (1,6) to (3,1);
		\draw (2,5) -- (4,3) -- (5,2);
		\draw [out=-15, in=90] (2,5) to (5,2);
		\draw [out=0, in=135] (2,5) to (6,4);
		\draw [out=270, in=135] (2,5) to (3,1);
		
		\plotperm{6,5,1,3,2,4}
		\draw[darkgray,thick,rounded corners=0.01,line cap=round]
			(0.5,0.5) rectangle (6.5,6.5);
		
		
		\draw[->] (7,3.5) to (8,3.5);
		
	\end{tikzpicture}
	\hspace{1mm}
	\begin{tikzpicture}[scale=0.5]
		\draw [draw = none, fill=gray!20] (2,0.5) rectangle (5,6.5);
		\draw [color = gray, dotted, thick, line cap=round] (2,0.5)--(2,6.5);
		\draw [color = gray, dotted, thick, line cap=round] (5,0.5)--(5,6.5);
		\draw [color = gray, dotted, thick, line cap=round] (0.5,2)--(6.5,2);
		\draw [color = gray, dotted, thick, line cap=round] (0.5,5)--(6.5,5);
		\draw [out=0, in=90] (1,6) to (6,4);
		\draw [out=-45, in=180] (1,6) to (5,5);
		\draw [out=-60, in=165] (1,6) to (4,3);
		\draw [out=-75, in=105] (1,6) to (2,2);
		\draw [out=-90, in=180] (1,6) to (3,1);
		\draw (2,2) -- (3,1);
		\draw (5,5) -- (6,4);
		
		\plotperm{6,2,1,3,5,4}
		\draw[darkgray,thick,rounded corners=0.01,line cap=round]
			(0.5,0.5) rectangle (6.5,6.5);
		
	\end{tikzpicture}
\end{center}
\end{footnotesize}
\caption{The inversion graph of $651324$ (left) and its image under left multiplication by~$t_{25}$~(right).}
\label{fig:inversion-graph-transposition}
\end{figure}
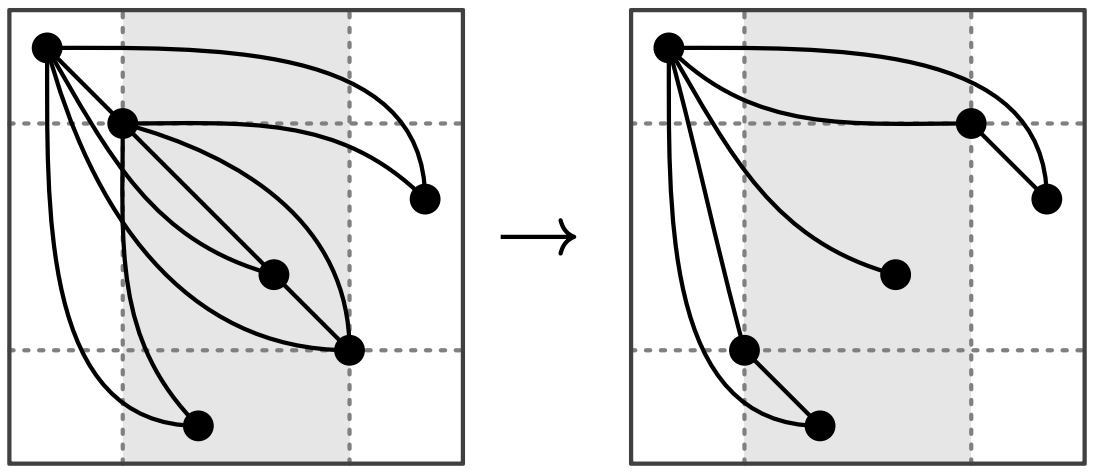

Throughout, reflections act on the left, so $t_{ab} \circ \pi$ is obtained from $\pi$ by swapping the values $a$ and $b$ in its one-line notation. We draw $G_\pi$ by placing the vertex $\pi(i)$ at the point $(i, \pi(i))$, so that two vertices are adjacent precisely when one lies southeast of the other; this embedding is a convenience for drawing, and inversion graphs do not come equipped with it. Figure~\ref{fig:inversion-graph-transposition} shows this for $\pi = 651324$ and $t_{25} \circ \pi$, with the shaded region marking the entries lying horizontally between the values $5$ and~$2$ in the plot of~$\pi$.

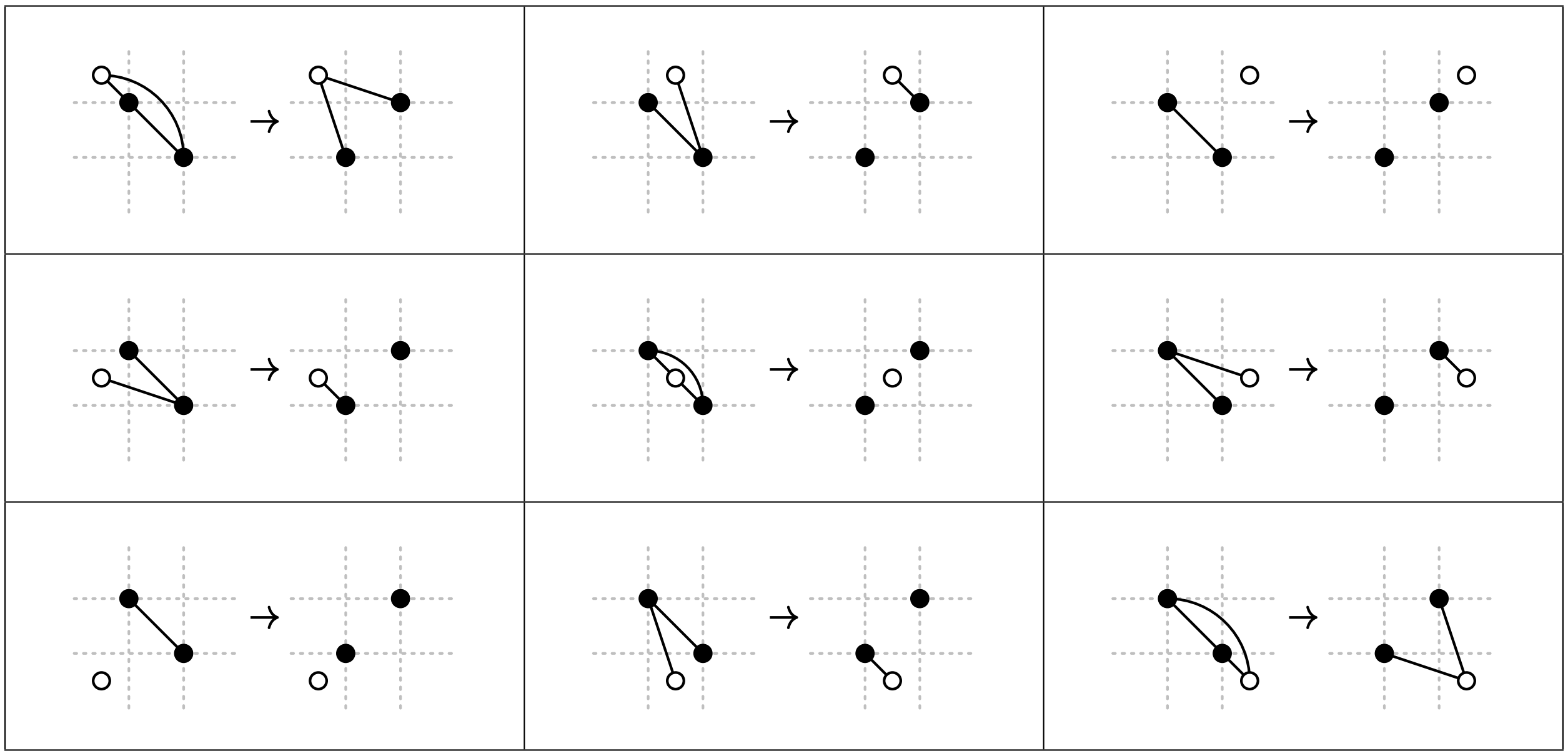
\begin{figure}[thbp]
\providecommand{\reflframe}[1]{%
  \begin{tikzpicture}[
    scale=0.58,
    baseline=(current bounding box.center)
  ]
    \draw[color=lightgray,dotted,thick,line cap=round]
      (1,0)--(1,3)
      (2,0)--(2,3)
      (0,1)--(3,1)
      (0,2)--(3,2);
    #1%
  \end{tikzpicture}%
}

\providecommand{\reflcase}[2]{%
  \raisebox{0pt}[\dimexpr\height+12pt\relax][\dimexpr\depth+12pt\relax]{%
    \mbox{%
      \reflframe{%
        \draw (1,2)--(2,1);
        #1
        \node[draw] at (1,2) [solid] {};
        \node[draw] at (2,1) [solid] {};
      }%
      \hspace{1mm}%
      \tikz[baseline=-0.5ex]{%
        \draw[->] (0,0)--(0.3,0);
      }%
      \hspace{1mm}%
      \reflframe{%
        #2
        \node[draw] at (1,1) [solid] {};
        \node[draw] at (2,2) [solid] {};
      }%
    }%
  }%
}

\renewcommand{\tabularxcolumn}[1]{%
  >{\centering\arraybackslash}m{#1}%
}

\begingroup

\setlength{\tabcolsep}{3pt}

\setlength{\arrayrulewidth}{0.4pt}

\noindent
\begin{tabularx}{\textwidth}{|X|X|X|}
\hline

\reflcase{
  \draw (0.5,2.5)--(1,2)
        (0.5,2.5) to[out=0,in=90] (2,1);
  \node[draw] at (0.5,2.5) [hollow] {};
}{
  \draw (2,2)--(0.5,2.5)--(1,1);
  \node[draw] at (0.5,2.5) [hollow] {};
}
&
\reflcase{
  \draw (1.5,2.5)--(2,1);
  \node[draw] at (1.5,2.5) [hollow] {};
}{
  \draw (1.5,2.5)--(2,2);
  \node[draw] at (1.5,2.5) [hollow] {};
}
&
\reflcase{
  \node[draw] at (2.5,2.5) [hollow] {};
}{
  \node[draw] at (2.5,2.5) [hollow] {};
}
\\
\hline

\reflcase{
  \draw (0.5,1.5)--(2,1);
  \node[draw] at (0.5,1.5) [hollow] {};
}{
  \draw (0.5,1.5)--(1,1);
  \node[draw] at (0.5,1.5) [hollow] {};
}
&
\reflcase{
  \draw (1,2) to[out=0,in=90] (2,1);
  \node[draw] at (1.5,1.5) [hollow] {};
}{
  \node[draw] at (1.5,1.5) [hollow] {};
}
&
\reflcase{
  \draw (2.5,1.5)--(1,2);
  \node[draw] at (2.5,1.5) [hollow] {};
}{
  \draw (2.5,1.5)--(2,2);
  \node[draw] at (2.5,1.5) [hollow] {};
}
\\
\hline

\reflcase{
  \node[draw] at (0.5,0.5) [hollow] {};
}{
  \node[draw] at (0.5,0.5) [hollow] {};
}
&
\reflcase{
  \draw (1.5,0.5)--(1,2);
  \node[draw] at (1.5,0.5) [hollow] {};
}{
  \draw (1.5,0.5)--(1,1);
  \node[draw] at (1.5,0.5) [hollow] {};
}
&
\reflcase{
  \draw (2.5,0.5)--(2,1)
        (1,2) to[out=0,in=90] (2.5,0.5);
  \node[draw] at (2.5,0.5) [hollow] {};
}{
  \draw (1,1)--(2.5,0.5)--(2,2);
  \node[draw] at (2.5,0.5) [hollow] {};
}
\\
\hline

\end{tabularx}

\endgroup
\caption{The effect of a reflection on an embedded inversion graph. The solid points represent an inversion being reflected; each panel shows a different position for a third entry, which is shown as a hollow point.}
\label{fig:inversion-graph-cases}
\end{figure}

Suppose that $a < b$ form an inversion of~$\pi$, so that $b$ appears to the left of~$a$. Passing from $\pi$ to $t_{ab} \circ \pi$ deletes the edge $ab$. A third entry~$c$ that does not lie horizontally between the swapped pair stays on the same side of both, so its adjacencies to $a$ and~$b$ are unaffected; if~$c$ lies horizontally between them, it changes sides with respect to both, and both adjacencies are toggled. Figure~\ref{fig:inversion-graph-cases} shows the nine positions available to a third entry and the effect of the swap in each.

One feature of the intermediate entries matters later: such an entry fails to be adjacent to~$b$ only if it is larger than~$b$, and fails to be adjacent to~$a$ only if it is smaller than~$a$; since~${a < b}$, it cannot do both, so it is adjacent to at least one of the pair.

The effect of the reflection is thus confined to the swapped pair and the entries horizontally between them: the edge joining the pair is deleted, and every adjacency between the pair and an intermediate entry is toggled. In lifting this to an arbitrary graph, we let any set of vertices play the intermediate role, subject only to the adjacency property just observed; this is the set~$X$ in the definition below.

All graphs in this paper are finite and simple. For a vertex $v$ of a graph~$G$, the \emph{neighborhood} $N_G(v)$ is the set of vertices adjacent to~$v$, and the \emph{closed neighborhood} is $N_G[v] = N_G(v) \cup \{v\}$. We drop the subscript when the graph is clear from context. Given two distinct vertices of~$G$, to \emph{toggle} their adjacency is to delete the edge joining them if it is present, and to add that edge otherwise.

\begin{definition}
\label{def:edge-reflection}
Let $G$ be a graph, let $u$ and~$v$ be distinct vertices of~$G$, and let $X$ be a set of vertices of~$G$ with $\{u, v\} \subseteq X$. The \emph{edge reflection} $t_{uv}^X$ transforms~$G$ into the graph $t_{uv}^X \circ G$ obtained by toggling the adjacency between $u$ and~$v$, together with the adjacencies of each $x \in X \setminus \{u, v\}$ to $u$ and~$v$.
\end{definition}

\begin{figure}[hbtp]
\begin{center}
\begin{tikzpicture}[scale=0.5,baseline=(current bounding box.center)]
	\plotpoints[w]{1/0,3/0,5/0,9/0}
	\plotpoints[v]{4/3,6/3}
	\labelnode[-90]{w1}{$w_1$}
	\labelnode[-90]{w2}{$w_2$}
	\labelnode[-90]{w3}{$w_3$}
    \node at (7,0) {$\dots$};
	\labelnode[-90]{w4}{\phantom{ww}$w_{|X|-2}$}
	\labelnode[135]{v1}{$u$}
	\labelnode[45]{v2}{$v$}
	\drawjoin{w1,w2,w3,w4,v1}{v1,v2}
\end{tikzpicture}
\end{center}
\caption{The pairs whose adjacencies are toggled by the edge reflection $t_{uv}^X$.}
\label{fig:nested-triangle}
\end{figure}
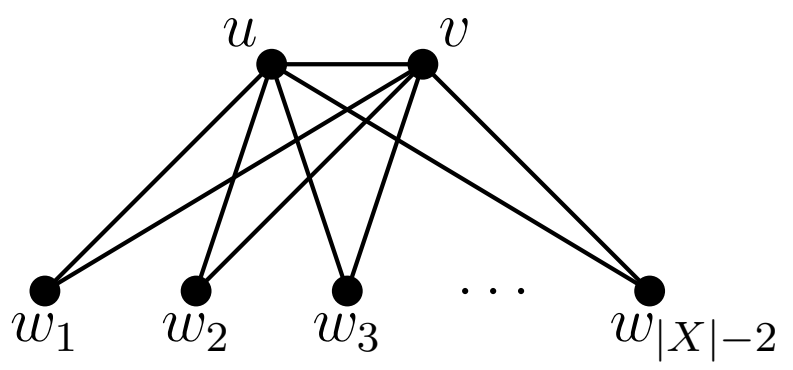

The pairs whose adjacencies are toggled by $t_{uv}^X$ are the edges of the graph $K_2 \vee \overline{K}_{|X|-2}$ on the vertex set~$X$, where $H_1 \vee H_2$ denotes the \emph{join}, obtained from disjoint copies of $H_1$ and~$H_2$ by adding every edge between them; here $u$ and~$v$ form the $K_2$, as in Figure~\ref{fig:nested-triangle}. Thus, the edge set of $t_{uv}^X \circ G$ is the symmetric difference of the edge sets of~$G$ and this join, and in particular $t_{uv}^X$ is an involution.

We apply edge reflections in one direction only.

\begin{definition}
\label{def:admissible}
The edge reflection $t_{uv}^X$ is \emph{admissible} at a graph~$G$ if $uv \in E(G)$ and $X \subseteq N(u) \cup N(v)$.
\end{definition}

Both conditions are satisfied by a permutation reflection applied to an inversion: the swapped pair forms an edge, and each intermediate entry is adjacent to at least one of them. Admissibility is a property of the reflection together with~$G$, not of the reflection alone. An admissible reflection deletes~$uv$, so if $t_{uv}^X$ is admissible at~$G$, it is not admissible at~$t_{uv}^X \circ G$.

Repeatedly applying admissible edge reflections brings a graph closer to the edgeless graph, as the following calculation shows. An admissible reflection $t_{uv}^X$ deletes $uv$, and each $x \in X \setminus \{u, v\}$ is adjacent to $u$ or~$v$; if $x$ is adjacent to both, it loses both edges, and otherwise it loses one edge and gains the other.

\begin{observation}
\label{obs:edge-loss}
If the edge reflection $t_{uv}^X$ is admissible at~$G$, then
\[
	|E(t_{uv}^X \circ G)| = |E(G)| - 2|X \cap N(u) \cap N(v)| - 1.
\]
\end{observation}

In particular, an admissible edge reflection strictly decreases the number of edges, and it does so by an odd number, a fact we return to in Section~\ref{sec:cyclomatic}.

We can now define the central quantity of this paper.

\begin{definition}
\label{def:edge-reflection-length}
A sequence of edge reflections \emph{empties} a graph if it transforms that graph into an edgeless graph and each reflection in the sequence is admissible at the graph to which it is applied. The \emph{edge reflection length} of~$G$, denoted by~$\graphrefl(G)$, counts the fewest edge reflections needed to empty~$G$.
\end{definition}

Because we require edge reflections to be admissible, the number of edges strictly decreases at each step of a sequence that empties a graph, and no step can undo an admissible reflection. Dropping this requirement can lower the number of reflections needed to empty a graph, in contrast to the situation for permutations, where sorting by reflections never benefits from a step that increases the number of inversions. Section~\ref{sec:further-questions} considers when this happens.

Every vertex whose adjacencies are toggled by an admissible reflection $t_{uv}^X$ lies in~$X$, and is therefore adjacent to $u$ or~$v$; thus the reflection acts within the component containing~$uv$, and in particular never joins two components. It follows that each reflection in a sequence that empties~$G$ acts within a single component of~$G$, so edge reflection length is additive.

\begin{observation}
\label{obs:additivity}
If $G$ has connected components $G_1, \ldots, G_k$, then
\[
	\graphrefl(G) = \graphrefl(G_1) + \cdots + \graphrefl(G_k).
\]
\end{observation}

Two upper bounds on $\graphrefl(G)$ follow from the definition, one counting edges and one counting vertices. Sections~\ref{sec:cyclomatic} and~\ref{sec:vertex-bound} show that each holds with equality precisely when $G$ is a forest.

\begin{proposition}
\label{prop:upper-bounds}
For every graph~$G$, we have
\begin{enumerate}[label=\textup{(\alph*)}]
	\item
	\label{prop:upper-bounds:edge}
	$\graphrefl(G) \le |E(G)|$, and
	\item
	\label{prop:upper-bounds:vertex}
	$\graphrefl(G) \le |V(G)| - \kappa(G)$.
\end{enumerate}
\end{proposition}

\begin{proof}
For~\ref{prop:upper-bounds:edge}, let $uv$ be any edge of~$G$. The reflection $t_{uv}^{\{u, v\}}$ is admissible at~$G$ and deletes $uv$ and nothing else, so~$G$ can be emptied one edge at a time.

For~\ref{prop:upper-bounds:vertex}, let $u$ be any nonisolated vertex of~$G$ and let $v$ be a neighbor of~$u$. The reflection~$t_{uv}^{N[u]}$ is admissible at~$G$ and isolates~$u$. A connected graph on~$k$ vertices can therefore be emptied by isolating its vertices one at a time, which takes at most $k - 1$ reflections (the last vertex is isolated automatically once all others are). Summing over the components of~$G$ using the additivity of Observation~\ref{obs:additivity} gives $\graphrefl(G) \le |V(G)| - \kappa(G)$.
\end{proof}

Finally, we record that edge reflections generalize permutation reflections, in the sense that reflecting an inversion of~$\pi$ and then passing to the inversion graph has the same effect as applying the corresponding edge reflection to~$G_\pi$.

\begin{observation}
\label{obs:compatibility}
Let $a$ and~$b$ form an inversion of a permutation~$\pi$, and let $X$ consist of $a$, $b$, and the entries lying horizontally between them in the plot of~$\pi$. Then $t_{ab}^X$ is admissible at~$G_\pi$ and
\[
	G_{t_{ab} \circ \pi} = t_{ab}^X \circ G_\pi.
\]
\end{observation}

Edge reflections are more general than permutation reflections in two respects. First, for an inversion $ab$ of~$\pi$, the permutation reflection $t_{ab}$ determines a single set~$X$ from the plot of~$\pi$, while $t_{ab}^X$ is admissible at~$G_\pi$ for every set~$X$ with $\{a, b\} \subseteq X \subseteq N(a) \cup N(b)$. Second, edge reflections apply to every graph, not only to inversion graphs. Section~\ref{sec:permutations} discusses how much of the permutation setting survives this additional freedom.

\section{The relation to sorting permutations by reflections}
\label{sec:sorting}

For a permutation~$\pi$, it may seem obvious that $\graphrefl(G_\pi) \le \ell'(\pi)$, especially in light of Observation~\ref{obs:compatibility}, but this does not follow immediately. The difficulty is that Observation~\ref{obs:compatibility} concerns only reflections applied to inversions of~$\pi$, each of which corresponds to an admissible edge reflection of~$G_\pi$; a reflection applied to a noninversion has no such counterpart, and, a priori, nothing prevents a shortest sorting sequence from using one. In this short section we show that such reflections are never needed: every permutation has a shortest sorting sequence that reflects only inversions.

This fact is the symmetric-group case of a theorem of Dyer~\cite{dyer:on-minimal-leng:}: in any Coxeter group, the reflection length of an element is achieved by some sequence of reflections each decreasing Coxeter length. In the symmetric group, decreasing Coxeter length means reflecting an inversion, so Dyer's theorem yields the sequence we need. Rather than invoke it, we give a short direct argument, and then use Observation~\ref{obs:compatibility} to convert the resulting sequence to edge reflections. Admissibility is the edge-reflection counterpart of this Coxeter-length decrease: it requires the reflected pair to be an edge, so every edge reflection we allow decreases the number of edges.

\begin{theorem}
\label{thm:graph-length-at-most-reflection-length}
For every permutation~$\pi$, we have $\graphrefl(G_\pi) \le \ell'(\pi)$.
\end{theorem}

\begin{proof}
We first recall why $\ell'(\pi) = n - c(\pi)$. Writing each $m$-cycle of~$\pi$ as a product of $m - 1$ transpositions expresses~$\pi$ using $n - c(\pi)$ of them, so $\ell'(\pi) \le n - c(\pi)$. Conversely, a transposition changes the number of cycles by one, so passing from the identity, with~$n$ cycles, to~$\pi$, with $c(\pi)$, takes at least $n - c(\pi)$ transpositions. Moreover, reflecting two values in a common cycle of~$\pi$ splits that cycle in two, so it decreases~$\ell'(\pi)$ by one.

We show, applying induction on $\ell'(\pi)$, that some shortest sequence sorting~$\pi$ reflects only inversions. If $\ell'(\pi) = 0$, then~$\pi$ is the identity and the empty sequence sorts it. Otherwise, $\pi$ has a cycle on at least two values. Let $a$ be the least value in such a cycle and put $b = \pi(a)$. Then $b > a$, and $b$ occupies position~$a$ while $a$ occupies position $\pi^{-1}(a) > a$, so $b$ appears to the left of~$a$. Thus, $ab$ is an inversion of~$\pi$ with $a$ and~$b$ in a common cycle. Reflecting such a pair decreases the reflection length, so $\ell'(t_{ab} \circ \pi) = \ell'(\pi) - 1$, and by induction some shortest sequence sorts $t_{ab} \circ \pi$ using only inversions. Since $ab$ is an inversion of~$\pi$, prepending $t_{ab}$ gives a shortest sequence sorting~$\pi$, again using only inversions.

Fix such a sequence, and let $\pi = \pi_0$, $\pi_1$, $\ldots$, $\pi_m$ be the permutations it passes through while sorting~$\pi$, where $m = \ell'(\pi)$ and $\pi_m$ is the identity. For each $i \in [m]$, the $i$th step reflects an inversion of~$\pi_{i-1}$, so by Observation~\ref{obs:compatibility} it corresponds to an admissible edge reflection sending $G_{\pi_{i-1}}$ to $G_{\pi_i}$. Since $G_{\pi_m}$ is edgeless, these $m$ edge reflections empty~$G_\pi$, and so $\graphrefl(G_\pi) \le m = \ell'(\pi)$.
\end{proof}

The inequality of Theorem~\ref{thm:graph-length-at-most-reflection-length} is often strict, and Section~\ref{sec:permutations} discusses two reasons that this can happen.

\section{Cyclomatic bounds}
\label{sec:cyclomatic}

A \emph{spanning forest} of a graph~$G$ is a spanning subgraph of~$G$ that is a forest, and a \emph{maximal spanning forest} is a spanning forest to which no edge of~$G$ can be added without creating a cycle. A maximal spanning forest therefore consists of a spanning tree of each connected component of~$G$, and so has $|V(G)| - \kappa(G)$ edges. The \emph{cyclomatic number} of~$G$ is
\[
	\beta(G) = |E(G)| - |V(G)| + \kappa(G),
\]
the number of edges that must be deleted to obtain a maximal spanning forest.

This section bounds the difference $|E(G)| - \graphrefl(G)$, which counts the reflections saved over emptying~$G$ one edge at a time, above and below in terms of $\beta(G)$. Our first lower bound follows simply by substitution into Proposition~\ref{prop:upper-bounds}\ref{prop:upper-bounds:vertex}:
\[
	|E(G)| - \graphrefl(G) \ge |E(G)| - |V(G)| + \kappa(G) = \beta(G).
\]
The upper bound $|E(G)| - \graphrefl(G) \le 2\beta(G)$ requires a bit more work. We start with the following result, which bounds the change in cyclomatic number under a single admissible edge reflection.

\begin{proposition}
\label{prop:one-step-cyclomatic}
If the edge reflection $t_{uv}^X$ is admissible at~$G$, then
\[
	\beta(t_{uv}^X \circ G) \le \beta(G) - |X\cap N(u)\cap N(v)|.
\]
\end{proposition}

\begin{proof}
Set $r = |X \cap N_G(u) \cap N_G(v)|$, so Observation~\ref{obs:edge-loss} states that
\[
	|E(G)| = |E(t_{uv}^X \circ G)| + 1 + 2r.
\]
Further take $F$ to be a maximal spanning forest of $t_{uv}^X \circ G$, so
\[
	\beta(t_{uv}^X \circ G) = |E(t_{uv}^X \circ G)| - |E(F)|,
\]
by definition. We construct a graph $F'$ with $|E(F')| = |E(F)| + 1 + r$ that contains a maximal spanning forest of~$G$. This suffices, because then
\begin{align*}
	\beta(G)
	&\ge
	|E(G)| - |E(F')|
	\\
	&=
	\left( |E(t_{uv}^X \circ G)| + 1 + 2r\right) - \left(|E(F)| + 1 + r\right)
	\\
	&=
	\beta(t_{uv}^X \circ G) + r.
\end{align*}

To construct $F'$ from $F$, we first add the edge $uv$. Each $x \in X \setminus \{u, v\}$ is adjacent in~$G$ to at least one of $u$ and~$v$, and $t_{uv}^X$ toggles both adjacencies, so $x$ is adjacent in $t_{uv}^X \circ G$, and hence in~$F$, to at most one of them. Then, for each $x\in X\setminus\{u,v\}$ that is adjacent in~$F$ to one of $u$ and~$v$, we toggle its adjacencies to $u$ and~$v$, deleting its edge to one of them and adding an edge to the other. Finally, for each of the $r$ vertices $x\in X\cap N_G(u)\cap N_G(v)$, we add the edge $xu$. Each edge we add joins a pair whose adjacency $t_{uv}^X$ toggles and that is not adjacent in $t_{uv}^X \circ G$, so it is an edge of~$G$ not already in~$F'$; each edge of~$F$ that we keep joins a pair whose adjacency $t_{uv}^X$ does not toggle, so it too is an edge of~$G$. Thus, $F'$ is a spanning subgraph of~$G$ with $|E(F')|=|E(F)|+r+1$, as claimed.

It remains to show that $F'$ contains a maximal spanning forest of $G$. Because $F'$ is a spanning subgraph of $G$, any maximal spanning forest of $F'$ is a spanning forest of $G$, and it is maximal in $G$ precisely when no edge of $G$ joins two components of $F'$. First, note that every edge of $F$ is either retained in $F'$ or replaced by a two-edge path passing through $uv$. Thus any two vertices in a common component of $F$, and hence in a common component of $t_{uv}^X \circ G$, lie in a common component of $F'$. Next, $u$ and $v$ are joined by the edge $uv$, each $x \in X \cap N_G(u) \cap N_G(v)$ is joined to $u$ by an added edge, and every other $x \in X \setminus \{u, v\}$ is adjacent in $t_{uv}^X \circ G$ to $u$ or~$v$; hence all vertices of $X$ lie in a single component of $F'$. Since every edge of $G$ either belongs to $t_{uv}^X \circ G$ or has both endpoints in $X$, no edge of $G$ joins two components of $F'$, as desired.
\end{proof}

Summing this bound along a sequence of edge reflections gives an upper bound on $|E(G)| - \graphrefl(G)$.

\begin{theorem}
\label{thm:sandwich}
For every graph~$G$, $|E(G)| - \graphrefl(G)$ is even, and
\[
	\beta(G) \le |E(G)| - \graphrefl(G) \le 2\beta(G).
\]
\end{theorem}

\begin{proof}
Set $n=|V(G)|$, $k=\graphrefl(G)$, and let $G = G_0, G_1, \ldots, G_k = \overline{K}_n$ be the graphs produced by a sequence of $k$ admissible edge reflections emptying $G$. For each~$i$, let $r_i$ be the number of common neighbors of the pair reflected at step~$i$ that lie in the set~$X$ of that reflection. By Observation~\ref{obs:edge-loss}, the $i$th reflection decreases the number of edges by $1 + 2r_i$, and since $G_k$ is edgeless,
\[
	|E(G)| = \sum_{i=1}^k (1 + 2r_i) = k + 2\sum_{i=1}^k r_i = \graphrefl(G) + 2\sum_{i=1}^k r_i.
\]
Thus $|E(G)| - \graphrefl(G) = 2\sum r_i$, which is even, proving the first claim of the theorem.

The lower bound on $|E(G)| - \graphrefl(G)$ has been established above, so it now suffices to prove the upper bound. By Proposition~\ref{prop:one-step-cyclomatic}, we have $\beta(G_{i-1})-\beta(G_i)\ge r_i$ for every $i$. Summing this inequality over $1\le i\le k$, the left-hand side telescopes, and we obtain
\[
	\beta(G_0) - \beta(G_k)
	=
	\beta(G) - \beta(\overline{K}_n)
	=
	\beta(G)
	\ge
	\sum_{i=1}^k r_i.
\]
Since $|E(G)| - \graphrefl(G) = 2\sum r_i$, this establishes that $|E(G)| - \graphrefl(G) \le 2\beta(G)$, completing the proof.
\end{proof}

Several consequences follow. First, the difference $|E(G)| - \graphrefl(G)$ is even, so it is never equal to~$1$. The same holds for permutations: the sign of~$\pi$ is both $(-1)^{\ell(\pi)}$ and $(-1)^{\ell'(\pi)}$, so $\ell(\pi) - \ell'(\pi)$ is even and never equal to~$1$.

\begin{corollary}
\label{cor:no-gap-one}
For every graph~$G$, we have $|E(G)| - \graphrefl(G) \neq 1$.
\end{corollary}

When $\beta(G) = 1$, the parity constraint also determines $\graphrefl(G)$: Theorem~\ref{thm:sandwich} gives $1\le |E(G)|-\graphrefl(G)\le 2$, and since this difference is even, it must equal~$2$.

\begin{corollary}
\label{cor:beta-one}
For every graph~$G$ with $\beta(G) = 1$, we have $|E(G)| - \graphrefl(G) = 2$. In particular, $\graphrefl(C_n) = n - 2$ for every $n \ge 3$.
\end{corollary}

The two bounds of Theorem~\ref{thm:sandwich} agree when $\beta(G) = 0$, identifying the graphs that cannot be emptied faster than one edge at a time.

\begin{corollary}
\label{cor:forest-edge-bound}
A graph~$G$ satisfies $\graphrefl(G) = |E(G)|$ if and only if~$G$ is a forest.
\end{corollary}

For forests, we can say more.

\begin{proposition}
\label{prop:forest-one-reflection}
Applying an admissible edge reflection to a forest yields a forest with one fewer edge.
\end{proposition}

\begin{proof}
Let $G$ be a forest and let $t_{uv}^X$ be admissible at~$G$, and put $H = t_{uv}^X \circ G$. Since $uv$ is an edge of~$G$ and $G$ has no triangle, the vertices $u$ and~$v$ have no common neighbor, so Observation~\ref{obs:edge-loss} gives $|E(H)| = |E(G)| - 1$. Proposition~\ref{prop:one-step-cyclomatic} gives $\beta(H) \le \beta(G) = 0$, so $H$ is also a forest.
\end{proof}

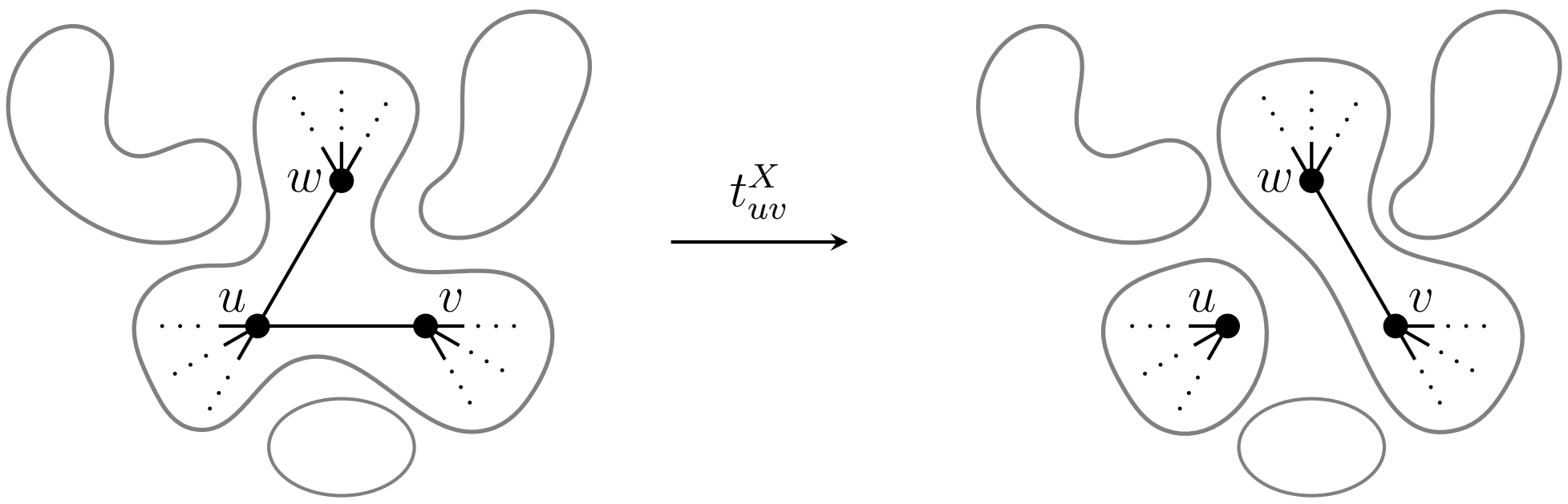
\begin{figure}[t]
\begin{center}

\newcommand{\attachthingy}[2]{%
	\draw
		(#1) -- ++(#2:0.8)
		node[pos=2.25,sloped,transform shape,scale=2] {$\dots$};
}
\newcommand{\attachthingies}[2]{
	\foreach \i in {#2} {
		\attachthingy{#1}{\i}
	}
}

\begin{tikzpicture}[>=stealth]

\node (T1) at (-6, 0) {
	\begin{tikzpicture}[scale=0.4,baseline=(current bounding box.center)]
		\node (u) at (210:2) [solid] {};
		\node (v) at (330:2) [solid] {};
		\node (w) at (90:2) [solid] {};
		\draw (w)--(u)--(v);
		\labelnode[180]{w}{$w$}
		\labelnode[135]{u}{$u$}
		\labelnode[45]{v}{$v$}
		\attachthingies{w}{60,90,120}
		\attachthingies{u}{180,210,240}
		\attachthingies{v}{300,330,0}
		\begin{scope}[on background layer]
		\draw[use Hobby shortcut, closed, line width=1pt, opacity=0.5]
		(40:1.0)..(70:4.2)..(90:4.5)..(110:4.3)..(130:2.5)..(150:1.8)..(170:2.0)..(175:2.8)..(190:4.3)..(210:4.5)..(225:4.4)..(240:2.0)..(315:4.4)..(330:4.6)..(345:4.5)..(360:3.6);
		\end{scope}
        \draw[use Hobby shortcut, closed, line width=1pt,opacity=0.5]
        (45:2.4)..(30:5.2)..(45:7.0)..(60:5.2)..(45:3.0);
        \draw[use Hobby shortcut, closed, line width=1pt,opacity=0.5]
        (135:3.0)..(165:6.0)..(135:7.0)..(150:5.2)..(135:4.0);
        \draw[opacity=0.5] (0,-3.5) ellipse (1.5cm and 1cm);
	\end{tikzpicture}
};

\node (T2) at (2, 0) {
	\begin{tikzpicture}[scale=0.4,baseline=(current bounding box.center)]
		\node (u) at (210:2) [solid] {};
		\node (v) at (330:2) [solid] {};
		\node (w) at (90:2) [solid] {};
		\draw (w)--(v);
		\labelnode[180]{w}{$w$}
		\labelnode[135]{u}{$u$}
		\labelnode[45]{v}{$v$}
		\attachthingies{w}{60,90,120}
		\attachthingies{u}{180,210,240}
		\attachthingies{v}{300,330,0}
		\begin{scope}[on background layer]
		\draw[use Hobby shortcut, closed, line width=1pt, opacity=0.5]
		(40:1.7)..(70:4.2)..(90:4.5)..(110:4.3)..(130:2.5)..(45:0.2)..(315:4.4)..(330:4.6)..(345:4.5)..(360:3.6);
		\draw[use Hobby shortcut, closed, line width=1pt, opacity=0.5]
		(170:2.0)..(175:2.8)..(190:4.3)..(210:4.5)..(225:4.4)..(240:2.0);
		\end{scope}
        \draw[use Hobby shortcut, closed, line width=1pt,opacity=0.5]
        (45:2.4)..(30:5.2)..(45:7.0)..(60:5.2)..(45:3.0);
        \draw[use Hobby shortcut, closed, line width=1pt,opacity=0.5]
        (135:3.0)..(165:6.0)..(135:7.0)..(150:5.2)..(135:4.0);
        \draw[opacity=0.5] (0,-3.5) ellipse (1.5cm and 1cm);
	\end{tikzpicture}
};

\draw[->, thick, shorten <=14pt, shorten >=14pt]
	(T1) -- node[above, yshift=4pt] {$t_{uv}^X$} (T2);

\end{tikzpicture}
\end{center}
\caption{An edge reflection $t_{uv}^X$ in a forest, with $X = \{u, v, w\}$. Enclosed regions indicate connected components (trees). The edges $uv$ and $uw$ are deleted and $vw$ is added, splitting the component and transferring the subtree containing~$w$ to~$v$'s component.}
\label{fig:reflection-on-forest}
\end{figure}

Figure~\ref{fig:reflection-on-forest} illustrates this. Because a reflection applied to a forest decreases the number of edges by exactly one, every sequence emptying a forest~$G$ has length $|E(G)|$. Forests are the only graphs whose emptying sequences all have the same length: if~$G$ contains a cycle, then emptying~$G$ one edge at a time takes~$|E(G)|$ reflections, which is more than~$\graphrefl(G)$ by Corollary~\ref{cor:forest-edge-bound}.

\section{Graphs meeting the vertex bound}
\label{sec:vertex-bound}

Corollary~\ref{cor:forest-edge-bound} shows that the edge bound $\graphrefl(G) \le |E(G)|$ holds with equality precisely when~$G$ is a forest. We now show that the same is true of the vertex bound $\graphrefl(G) \le |V(G)| - \kappa(G)$ of Proposition~\ref{prop:upper-bounds}\ref{prop:upper-bounds:vertex}.

For a forest, the edge and vertex bounds coincide, so forests attain equality in the vertex bound. It remains only to prove the converse, that a graph with a cycle does not.

\begin{theorem}
\label{thm:forest-iff-vertex-bound}
A graph~$G$ satisfies $\graphrefl(G) = |V(G)| - \kappa(G)$ if and only if~$G$ is a forest.
\end{theorem}

\begin{proof}
Both $\graphrefl(G)$ and $|V(G)| - \kappa(G)$ are additive over connected components, so it suffices to establish that a connected graph~$G$ on $n$ vertices containing a cycle satisfies $\graphrefl(G) \le n - 2$, which then applies to any component containing a cycle.

Isolating nonisolated vertices one at a time gives $\graphrefl(G) \le n - 1$, as in the proof of Proposition~\ref{prop:upper-bounds}\ref{prop:upper-bounds:vertex}. The same count applies to any sequence of $m$ reflections emptying~$G$ in which each reflection isolates at least one new vertex: an isolated vertex never lies in the set~$X$ of an admissible reflection and so stays isolated, and the last reflection isolates both endpoints of an edge; hence $m \le n - 1$. Every reflection in the procedure below isolates at least one new vertex. To show that it uses at most $n - 2$ reflections, we either find a reflection before the last that isolates at least two vertices or show that the last reflection isolates at least three.

Fix an induced cycle $C$ in~$G$, with vertices $c_1, c_2, \ldots, c_k$ in cyclic order, and label the remaining vertices $v_1, \ldots, v_{n-k}$ arbitrarily. Isolating the $v_i$ one by one to leave only the cycle and then invoking Corollary~\ref{cor:beta-one} is tempting but does not work: isolating $v_i$ along an edge to a cycle vertex can change edges among the cycle vertices and so disturb the cycle.

We instead begin by clearing the edges among the $v_i$. We process $v_1, \ldots, v_{n-k}$ in order. At step~$i$, if $v_i$ is adjacent to some~$v_j$ with~${j > i}$, we apply~$t_{v_i v_j}^{N[v_i]}$ to delete every edge at~$v_i$ and so isolate it; otherwise we do nothing.

By induction on~$h$, once we have processed~$v_h$, it has no edges to the other~$v_i$. Indeed, if $v_h$ was reflected, that reflection isolated it. If it was skipped, it had no neighbor~$v_j$ with $j > h$, or it would have been reflected, and by the induction hypothesis it had no neighbor~$v_i$ with $i < h$. Furthermore, no later step adds an edge between~$v_h$ and another $v_i$, since a step at~$v_{i'}$ with $i' > h$ toggles only adjacencies between vertices of $N[v_{i'}]$, and~$v_h$ does not lie in that set because it is not adjacent to $v_{i'}$.

After this pass, no edges remain among $v_1, \ldots, v_{n-k}$, every reflection performed has isolated at least one new vertex, and~$C$ is still an induced cycle, since every adjacency toggled in the pass involves one of the~$v_i$. If~$v_i$ still has edges, they all join it to cycle vertices; we call such~$v_i$ \emph{outer vertices}. Figure~\ref{fig:tree_iff_n-1:phase-1} shows the resulting graph, ignoring isolated vertices.

\begin{figure}[t]
\begin{center}
	\begin{tikzpicture}[scale=4.4, baseline=(current bounding box.center)]
    \makeninegon{}{}{}
	\draw [ultra thick] (v6)--(v7)--(v8)--(v9)--(v1)--(v2)--(v3)--(v4)--(v5)--(v6);

	\node[draw] (l1) at (-.55,.35) [hollow] {};
	\drawjoin{l1}{v1,v2,v4}
	\node[draw] (l2) at (-.65,.05) [hollow] {};
	\drawjoin{l2}{v2,v4}
	\node[draw] (l3) at (-.5,-.3) [hollow] {};
	\drawjoin{l3}{v3,v4,v5}

	\node[draw] (r1) at (.6,.35) [hollow] {};
	\drawjoin{r1}{v1,v7}
	\node[draw] (r2) at (.65,.11) [hollow] {};
	\drawjoin{r2}{v7,v8,v9}
	\node[draw] (r3) at (.55,-.3) [hollow] {};
	\drawjoin{r3}{v5,v7,v8}

	\node[draw] (c1) at (-.075,-.05) [hollow] {};
	\drawjoin{c1}{v1,v2,v3,v4,v5,v7}
	\node[draw] (c2) at (.055,.08) [hollow] {};
	\drawjoin{c2}{v1,v3,v6,v8}


	\end{tikzpicture}
\end{center}
\caption{The structure of the graph after the initial pass, ignoring isolated vertices: an induced cycle (solid dots) together with an independent set of outer vertices (hollow dots).}
\label{fig:tree_iff_n-1:phase-1}
\end{figure}
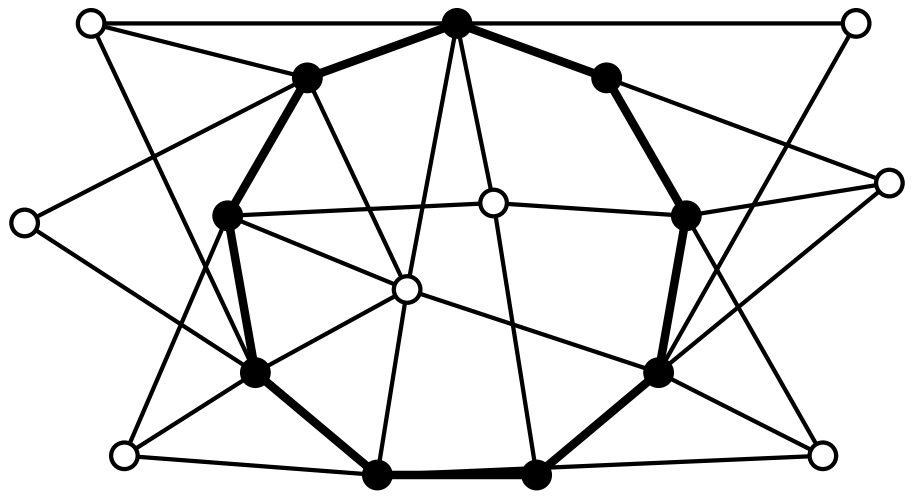

Next we shrink the cycle to a triangle. If~${k \ge 4}$, we process $c_1, \ldots, c_{k-3}$ in order, applying~$t_{c_i c_{i+1}}^{N[c_i]}$ at step~$i$. This isolates~$c_i$, and among $c_{i+1}, \ldots, c_k$ the only change is that the edge~$c_{i+1} c_k$ is added, since $c_k$ was adjacent to~$c_i$ but not to~$c_{i+1}$; so these vertices still induce a cycle. Outer vertices may become isolated along the way, which can only help. After step~${k - 3}$, the vertices $c_{k-2}, c_{k-1}, c_k$ induce a triangle, and we relabel them $p, q, r$.

If an outer vertex is adjacent to exactly one of $p, q, r$, we isolate it with a single reflection along that edge. After these reflections, each remaining outer vertex is adjacent to at least two of $p, q, r$. We group these vertices by their adjacencies to the triangle, setting
\begin{align*}
	P &= \bigl(N(q) \cap N(r)\bigr) \setminus N[p], \\
	Q &= \bigl(N(p) \cap N(r)\bigr) \setminus N[q], \\
	R &= \bigl(N(p) \cap N(q)\bigr) \setminus N[r], \\
	A &= N(p) \cap N(q) \cap N(r),
\end{align*}
as depicted on the left of Figure~\ref{fig:tree_iff_n-1:endgame}. The vertices $p, q, r$ and the sets $P, Q, R$ play symmetric roles, so we may relabel the triangle to ensure that $|P| \ge |Q| \ge |R|$; then $R$ is empty whenever $Q$ is, and $P$ is nonempty whenever $Q$ is.

\begin{figure}[t]
\begin{center}
\begin{tikzpicture}[>=stealth]

\newcommand{\treefigcommon}{
	\node[draw] (p) at (0,{sqrt(2)}) [solid] {};
	\node[draw] (q) at (-1,0) [solid] {};
	\node[draw] (r) at (1,0) [solid] {};
	\labelnode[45]{p}{$p$}
	\labelnode[135]{q}{$q$}
	\labelnode[45]{r}{$r$}

	\begin{scope}[shift={(0,-.5)}]
	\begin{scope}[on background layer]
	\draw[use Hobby shortcut, closed, line width=1pt, opacity=0.5, fill=lightgray]
	(0:.4)..(45:.3)..(135:.35)..(270:.32);
	\end{scope}
	\foreach \angle/\dist [count=\i] in {0/.1, 110/.22, 270/.2} {
		\node[draw] (P\i) at (\angle:\dist) [hollow] {};
	}
	\node at (345:.53) {$P$};
	\end{scope}

	\begin{scope}[shift={({.5+sqrt(2)/4},{3*sqrt(2)/4})}]
	\begin{scope}[on background layer]
	\draw[use Hobby shortcut, closed, line width=1pt, opacity=0.5, fill=lightgray]
	(0:.25)..(45:.3)..(150:.34)..(225:.28)..(280:.34)..(325:.27);
	\end{scope}
	\foreach \angle/\dist [count=\i] in {60/.18, 150/.19, 280/.2} {
		\node[draw] (Q\i) at (\angle:\dist) [hollow] {};
	}
	\node at (-50:.45) {$Q$};
	\end{scope}

	\begin{scope}[shift={({-.5-sqrt(2)/4},{3*sqrt(2)/4})}]
	\begin{scope}[on background layer]
	\draw[use Hobby shortcut, closed, line width=1pt, opacity=0.5, fill=lightgray]
	(0:.25)..(45:.3)..(90:.28)..(135:.31)..(270:.34)..(315:.27);
	\end{scope}
	\foreach \angle/\dist [count=\i] in {0/.1, 160/.18, 250/.22} {
		\node[draw] (R\i) at (\angle:\dist) [hollow] {};
	}
	\node at (220:.51) {$R$};
	\end{scope}

	\begin{scope}[shift={(0,{sqrt(2)/2})}]
	\begin{scope}[on background layer]
	\draw[use Hobby shortcut, closed, line width=1pt, opacity=0.5, fill=lightgray]
	(0:.25)..(45:.23)..(90:.28)..(135:.31)..(270:.36)..(315:.27);
	\end{scope}
	\foreach \angle/\dist [count=\i] in {0/.1, 160/.18, 250/.22} {
		\node[draw] (A\i) at (\angle:\dist) [hollow] {};
	}
	\node at (270:.47) {$A$};
	\end{scope}
}

\node (L) at (-8, 0) {
	\begin{tikzpicture}[scale=2, baseline=(current bounding box.center)]
	\treefigcommon 
	\draw [ultra thick] (r)--(p)--(q)--(r);
	\foreach \i in {1,2,3} {
		\draw (q)--(P\i)--(r);
		\draw (p)--(Q\i)--(r);
		\draw (p)--(R\i)--(q);
		\draw (p)--(A\i)--(q) (A\i)--(r);
	}
	\end{tikzpicture}
};

\node (R) at (0, 0) {
	\begin{tikzpicture}[scale=2, baseline=(current bounding box.center)]
	\treefigcommon 
	\draw [ultra thick] (r)--(p)--(q);
	\foreach \i in {1,2,3} {
		\draw (p)--(Q\i)--(r);
		\draw (p)--(R\i)--(q);
		\draw (p)--(A\i);
	}
	\end{tikzpicture}
};

\draw[->, thick, shorten <=14pt, shorten >=14pt]
	(L) -- node[above, yshift=4pt] {$t_{qr}^{\{q,r\} \cup P \cup A}$} (R);

\end{tikzpicture}
\end{center}
\caption{On the left, the situation after reducing to a triangle $\{p,q,r\}$ with outer vertices adjacent to at least two vertices of the triangle. On the right, the result of the first edge reflection in the case where $P,Q\neq\emptyset$.}
\label{fig:tree_iff_n-1:endgame}
\end{figure}
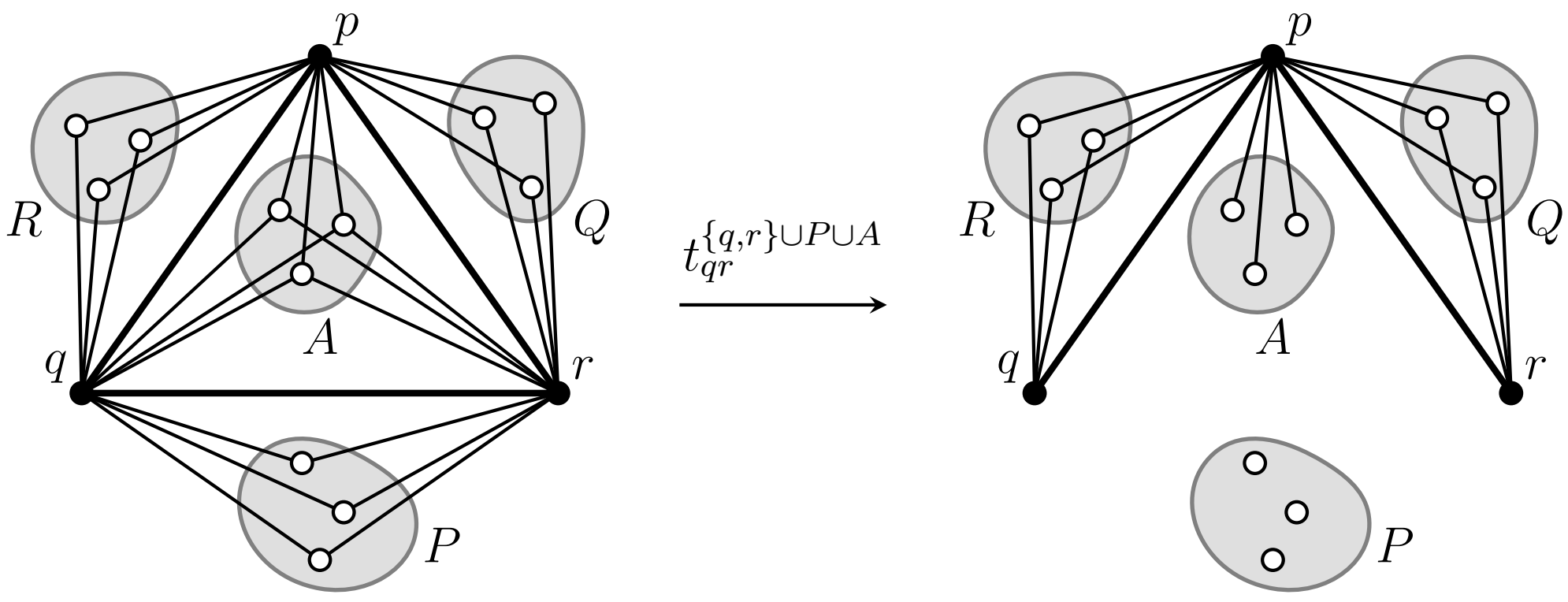

Every reflection so far has isolated at least one new vertex, so it remains to find the extra isolation promised at the start of the proof. We consider two cases. In each, after the reflections we describe, we finish by isolating the remaining nonisolated vertices one at a time; each of these reflections also isolates at least one new vertex.

First, suppose that $Q = \emptyset$, so that $R = \emptyset$ as well and the outer vertices are those of $P$ and~$A$. The reflection $t_{qr}^{\{p, q, r\} \cup P \cup A}$ deletes all three triangle edges and all edges from $q$ and~$r$ to $P$ and~$A$, isolating $q$, $r$, and every vertex of~$P$. It does not toggle the adjacencies of $p$ to the vertices of~$A$, so each vertex of~$A$ keeps its edge to~$p$ and no other. If $A \ne \emptyset$, some edge remains, so the reflection is not the last, and it isolates at least two vertices, $q$ and~$r$. If $A = \emptyset$, the reflection empties the graph, so it is the last, and it isolates at least three vertices, $p$, $q$, and~$r$.

Otherwise, $Q \neq \emptyset$, and thus $P \neq \emptyset$ as well. We apply $t_{qr}^{\{q, r\} \cup P \cup A}$, which deletes $qr$ and all edges from $q$ and~$r$ to $P$ and~$A$, isolating every vertex of~$P$ and giving the graph on the right of Figure~\ref{fig:tree_iff_n-1:endgame}. We then apply $t_{pr}^{\{p, r\} \cup Q}$, which deletes $pr$ and all edges from $p$ and~$r$ to~$Q$, isolating~$r$ and every vertex of~$Q$. Neither reflection toggles the adjacency of $p$ and~$q$, so the edge $pq$ remains. Thus, the second reflection is not the last, and it isolates at least two vertices, $r$ and the vertices of~$Q$.
\end{proof}

The proof of Theorem~\ref{thm:forest-iff-vertex-bound} gives more than the characterization itself: its procedure empties every connected graph~$H$ containing a cycle using at most $|V(H)| - 2$ reflections. More generally, for any graph~$G$, applying this to each component of~$G$ that contains a cycle, and Proposition~\ref{prop:upper-bounds}\ref{prop:upper-bounds:vertex} to the rest, gives the following.

\begin{corollary}
\label{cor:vertex-bound-improved}
For every graph~$G$, we have $\graphrefl(G) \le |V(G)| - \kappa(G) - q(G)$, where $q(G)$ is the number of connected components of~$G$ that contain a cycle.
\end{corollary}

Because $|E(G)| = |V(G)| - \kappa(G) + \beta(G)$, Corollary~\ref{cor:vertex-bound-improved} yields the following strengthening of Theorem~\ref{thm:sandwich}, identical to it except that the lower bound increases by~$q(G)$.

\begin{theorem}
\label{thm:sandwich-improved}
For every graph~$G$, $|E(G)| - \graphrefl(G)$ is even, and
\[
	\beta(G) + q(G) \le |E(G)| - \graphrefl(G) \le 2\beta(G),
\]
where $q(G)$ is the number of connected components of~$G$ that contain a cycle.
\end{theorem}

\begin{proof}
The upper bound and the parity statement are Theorem~\ref{thm:sandwich}. For the lower bound, Corollary~\ref{cor:vertex-bound-improved} gives
\[
	|E(G)| - \graphrefl(G) \ge |E(G)| - |V(G)| + \kappa(G) + q(G) = \beta(G) + q(G). \qedhere
\]
\end{proof}

The bounds are close for connected graphs with small cyclomatic numbers, and parity then determines the difference exactly. Applying this componentwise gives an exact value.

\begin{corollary}
\label{cor:exact-small-beta}
If every connected component of~$G$ has cyclomatic number at most~$2$, then $|E(G)| - \graphrefl(G) = 2\beta(G)$.
\end{corollary}

\begin{proof}
By Observation~\ref{obs:additivity}, both $|E(G)| - \graphrefl(G)$ and $2\beta(G)$ are sums of the corresponding quantities over the connected components of~$G$, so it suffices to consider a connected graph~$H$ with $\beta(H) \le 2$. If $\beta(H) = 0$, then $H$ is a tree and Corollary~\ref{cor:forest-edge-bound} gives $|E(H)| - \graphrefl(H) = 0$. Otherwise $H$ contains a cycle, so $q(H) = 1$ and Theorem~\ref{thm:sandwich-improved} gives
\[
	\beta(H) + 1 \le |E(H)| - \graphrefl(H) \le 2\beta(H).
\]
When $\beta(H) = 1$ these bounds are both~$2$. When $\beta(H) = 2$ they are $3$ and~$4$, and the difference is even by Theorem~\ref{thm:sandwich-improved}, hence equal to~$4$. In each case $|E(H)| - \graphrefl(H) = 2\beta(H)$.
\end{proof}

\section{Lower bounds and exact values}
\label{sec:lower-bounds}

The bounds of Sections~\ref{sec:cyclomatic} and~\ref{sec:vertex-bound} determine $\graphrefl(G)$ when $\beta(G)$ is small, but they say little about graphs with many independent cycles. This section develops two lower bounds in terms of standard graph parameters and then collects the exact values they yield.

We begin with the graphs that can be emptied in a single reflection. Suppose that $t_{uv}^X$ empties~$G$. Then every edge of~$G$ is among the pairs whose adjacencies it toggles, and every such pair must have been an edge of~$G$, since any pair that was not would be an edge afterwards; so each $x \in X \setminus \{u, v\}$ is adjacent to both $u$ and~$v$.

\begin{observation}
\label{obs:nested-triangles}
A graph~$G$ satisfies $\graphrefl(G) = 1$ if and only if $G$ consists of $K_2 \vee \overline{K}_m$ for some $m \ge 0$, together with isolated vertices.
\end{observation}

Our first lower bound counts the vertices that a sequence of reflections must touch. We include it, even though Proposition~\ref{prop:edge-domination-lower-bound} below is stronger, because its proof is a direct count, and the proof of Proposition~\ref{prop:edge-domination-lower-bound} refines the same idea. A \emph{vertex cover} of~$G$ is a set of vertices $S \subseteq V(G)$ such that every edge of~$G$ has at least one endpoint in~$S$, and $\tau(G)$ denotes the minimum size of a vertex cover.

\begin{proposition}
\label{prop:vertex-cover-lower-bound}
For every graph~$G$, we have $\graphrefl(G) \ge \lceil \tau(G)/2 \rceil$.
\end{proposition}

\begin{proof}
Suppose that $k$ edge reflections empty~$G$, and let $S$ be the set of vertices appearing as endpoints of the reflected pairs, so that $|S| \le 2k$. Every pair whose adjacency is toggled by a reflection has an endpoint in~$S$, so an edge of~$G$ with neither endpoint in~$S$ would remain an edge throughout the sequence. Hence $S$ is a vertex cover of~$G$, giving $\tau(G) \le 2k$.
\end{proof}

Our second lower bound counts edges rather than vertices. Its proof requires an argument about a single reflection, rather than a direct count. A set $D \subseteq E(G)$ is an \emph{edge dominating set} of~$G$ if every edge of~$G$ shares an endpoint with some edge of~$D$, and the \emph{edge domination number} $\gamma'(G)$ is the minimum size of an edge dominating set, which is~$0$ when $G$ is edgeless.

\begin{proposition}
\label{prop:edge-domination-lower-bound}
For every graph~$G$, we have $\graphrefl(G) \ge \gamma'(G)$.
\end{proposition}

\begin{proof}
It suffices to show that a single admissible edge reflection decreases the edge domination number by at most one, since applying this at each step of a sequence of $k$ reflections emptying~$G$ gives $\gamma'(G) \le k$.

Let $t_{uv}^X$ be admissible at~$G$, put $H = t_{uv}^X \circ G$, and let $D$ be a minimum edge dominating set of~$H$. We build from~$D$ an edge dominating set of~$G$ with at most one more edge. For each edge~${e \in D}$, define an edge~$e'$ of~$G$ as follows: if $e$ joins one of $u$ and~$v$ to some $w \in X \setminus \{u, v\}$, let $e'$ join the other of $u$ and~$v$ to~$w$, and otherwise put $e' = e$. Since $t_{uv}^X$ toggles the adjacencies of $w$ to $u$ and~$v$, and $w$ is adjacent in~$G$ to at least one of them, each such~$e'$ is an edge of~$G$.

We claim that $D' = \{e' : e \in D\} \cup \{uv\}$ is an edge dominating set of~$G$. Every edge of~$G$ incident with $u$ or~$v$ shares an endpoint with $uv \in D'$. Now let $f$ be an edge of~$G$ incident with neither~$u$ nor~$v$. Then~$f$ is an edge of~$H$, so it shares an endpoint with some $e \in D$. If $e' = e$, then $f$ shares an endpoint with $e' \in D'$. Otherwise $e$ and $e'$ have the same endpoint outside $\{u, v\}$ and differ only in whether the other endpoint is~$u$ or~$v$; since $f$ is incident with neither~$u$ nor~$v$, the endpoint it shares with~$e$ is the one that $e$ and $e'$ have in common, so $f$ shares an endpoint with $e'$ as well.

Thus $D'$ is an edge dominating set of~$G$ with $|D'| \le \gamma'(H) + 1$, giving $\gamma'(G) \le \gamma'(H) + 1$, which is the required bound.
\end{proof}

The endpoints of an edge dominating set of size~$k$ form a vertex cover of size at most~$2k$, so $\tau(G) \le 2\gamma'(G)$ and thus Proposition~\ref{prop:edge-domination-lower-bound} implies Proposition~\ref{prop:vertex-cover-lower-bound}. The improvement can be substantial. For the matching~$M_k$ consisting of $k$ disjoint edges, we have ${\tau(M_k) = \gamma'(M_k) = k}$, so the vertex cover bound gives only $\lceil k/2 \rceil$ while the edge domination bound gives~$k$, which is the exact value by Proposition~\ref{prop:upper-bounds}\ref{prop:upper-bounds:edge}.

Neither bound is sharp in general. The star $K_{1,m}$ has $\gamma'(K_{1,m}) = 1$ and $\tau(K_{1,m}) = 1$, so both bounds give~$1$, while $\graphrefl(K_{1,m}) = m$ by Corollary~\ref{cor:forest-edge-bound}.

We close by considering complete graphs and balanced complete bipartite graphs, where these bounds suffice to determine exact values of $\graphrefl$. For the complete graph, $\tau(K_n) = n - 1$, so Proposition~\ref{prop:vertex-cover-lower-bound} gives $\graphrefl(K_n) \ge \lfloor n/2 \rfloor$. This is the exact value: the reflection $t_{uv}^{N(u) \cup N(v)}$ isolates both $u$ and~$v$ at once, and repeating it empties $K_n$ in~$\lfloor n/2 \rfloor$ steps.

For the complete bipartite graph $K_{n,n}$, a perfect matching is an edge dominating set, and no smaller set dominates every edge, since a set of fewer than $n$ edges misses a vertex on each side, and the edge joining two such vertices is not dominated. Thus, $\gamma'(K_{n,n}) = n$, and Proposition~\ref{prop:edge-domination-lower-bound} gives $\graphrefl(K_{n,n}) \ge n$. This is the exact value, and the quickest route to the upper bound is through a permutation realizing the graph. The graph $K_{n,n}$ is the inversion graph of the permutation $(n+1)(n+2) \cdots (2n)\, 12 \cdots n$, which is the product of the $n$ disjoint transpositions $(1, n+1), \ldots, (n, 2n)$. Thus, we have
\[
	\graphrefl(K_{n,n})
	=
	\graphrefl(G_{(n+1) \cdots (2n)\, 1 \cdots n})
	\le
	\ell'((n+1) \cdots (2n)\, 1 \cdots n)
	=
	n.
\]

\section{Back to permutations}
\label{sec:permutations}

We now return to permutations and ask what the graph results say there. Since Coxeter length is the number of edges, we have $\ell(\pi) = |E(G_\pi)|$, and the quantity $|E(G)| - \graphrefl(G)$ that Theorem~\ref{thm:sandwich} bounds becomes $\ell(\pi) - \graphrefl(G_\pi)$. This sits alongside the permutation length gap $\ell(\pi) - \ell'(\pi)$, and the two are ordered by Theorem~\ref{thm:graph-length-at-most-reflection-length}, since $\graphrefl(G_\pi) \le \ell'(\pi)$ gives $\ell(\pi) - \ell'(\pi) \le \ell(\pi) - \graphrefl(G_\pi)$. Tenner and Vatter~\cite[Corollary 1.5]{TV} bound the length gap below by $\beta(G_\pi)$, and Theorem~\ref{thm:sandwich} bounds ours above by $2\beta(G_\pi)$.

\begin{theorem}
\label{thm:permutation-chain}
For every permutation~$\pi$,
\[
	\beta(G_\pi) \le \ell(\pi) - \ell'(\pi) \le \ell(\pi) - \graphrefl(G_\pi) \le 2\beta(G_\pi).
\]
\end{theorem}

\begin{proof}
The first inequality is due to Tenner and Vatter~\cite[Corollary 1.5]{TV}. The second is Theorem~\ref{thm:graph-length-at-most-reflection-length}, and the third is Theorem~\ref{thm:sandwich} applied to~$G_\pi$, using $\ell(\pi) = |E(G_\pi)|$.
\end{proof}

This chain of inequalities bounds how far the edge reflection length falls below the reflection length. Writing
\[
	\ell'(\pi) - \graphrefl(G_\pi) = \bigl(\ell(\pi) - \graphrefl(G_\pi)\bigr) - \bigl(\ell(\pi) - \ell'(\pi)\bigr),
\]
the middle inequality of the chain makes this difference nonnegative, and the outer two bound it above by $2\beta(G_\pi) - \beta(G_\pi) = \beta(G_\pi)$.

\begin{corollary}
\label{cor:reflection-gap}
For every permutation~$\pi$,
\[
	0 \le \ell'(\pi) - \graphrefl(G_\pi) \le \beta(G_\pi).
\]
\end{corollary}

Moreover, this difference is even. Both $\ell(\pi) - \graphrefl(G_\pi)$ and $\ell(\pi) - \ell'(\pi)$ are even, the first by Theorem~\ref{thm:sandwich} and the second as noted after that theorem, so their difference $\ell'(\pi) - \graphrefl(G_\pi)$ is even as well.

The upper bound can be sharpened, because Tenner and Vatter~\cite[Theorem 1.4]{TV} decompose the permutation length gap exactly rather than merely bounding it:
\[
	\ell(\pi) - \ell'(\pi) = \beta(G_\pi) + \bigl(c(\pi) - \kappa(\pi)\bigr),
\]
where $\kappa(\pi)$ is the number of \emph{sum components} of~$\pi$, the parts in its finest expression as a direct sum $\pi = \pi_1 \oplus \cdots \oplus \pi_{\kappa(\pi)}$, in which each~$\pi_i$ acts on a contiguous block of values above the previous block. The sum components of~$\pi$ correspond to the connected components of~$G_\pi$, so $\kappa(\pi) = \kappa(G_\pi)$ and the notation is consistent with our use of~$\kappa$ for graphs. The quantity $c(\pi) - \kappa(\pi)$ is the number of cycles of~$\pi$ in excess of its sum components, which Tenner and Vatter term the \emph{cycle surplus} of~$\pi$. Subtracting this identity from $\ell(\pi) - \graphrefl(G_\pi) \le 2\beta(G_\pi)$ gives the following, whose right-hand side is nonnegative because the cycle surplus is at most $\beta(G_\pi)$ by Tenner and Vatter~\cite[Theorem 1.3]{TV}.

\begin{corollary}
\label{cor:refined-gap}
For every permutation~$\pi$,
\[
	\ell'(\pi) - \graphrefl(G_\pi) \le \beta(G_\pi) - \bigl(c(\pi) - \kappa(\pi)\bigr).
\]
In particular, $\ell'(\pi)$ and $\graphrefl(G_\pi)$ agree whenever the cycle surplus of~$\pi$ equals $\beta(G_\pi)$.
\end{corollary}

The criterion is met whenever every connected component of~$G_\pi$ has cyclomatic number at most~$1$. To see this, note that the cycle surplus and $\beta(G_\pi)$ are both sums over the sum components of~$\pi$, which correspond to the connected components of~$G_\pi$; for the cycle surplus, this holds because each sum component is $\pi$-invariant, so every cycle of~$\pi$ lies in a single sum component. It therefore suffices to consider $\pi$ with $G_\pi$ connected and $\beta(G_\pi) \le 1$. If $\beta(G_\pi) = 0$, the cycle surplus is~$0$, matching $\beta(G_\pi)$. If $\beta(G_\pi) = 1$, the cycle surplus is at most~$1$, and it has the same parity as $\beta(G_\pi)$ because their sum $\ell(\pi) - \ell'(\pi)$ is even, so the surplus is~$1$, again matching $\beta(G_\pi)$.

\begin{corollary}
\label{cor:unicyclic-equality}
If every connected component of~$G_\pi$ has cyclomatic number at most~$1$, then $\ell'(\pi) = \graphrefl(G_\pi)$.
\end{corollary}

The case $\beta(G_\pi) = 0$ deserves emphasis, because there the chain of Theorem~\ref{thm:permutation-chain} collapses entirely.

\begin{corollary}
\label{cor:forest-all-equal}
If $G_\pi$ is a forest, then $\graphrefl(G_\pi) = \ell'(\pi) = \ell(\pi) = |E(G_\pi)|$.
\end{corollary}

It is folklore that $G_\pi$ is a forest if and only if~$\pi$ avoids $321$ and $3412$; see~\cite[Proposition 2.4]{TV} for a proof. Corollary~\ref{cor:forest-all-equal} therefore gives one direction of the theorem of Petersen and Tenner~\cite{petersen:the-depth-of-a-:} that $\ell(\pi) = \ell'(\pi)$ holds precisely on $\mathrm{Av}(321, 3412)$: if $\pi$ avoids the two patterns, then $G_\pi$ is a forest, so $\graphrefl(G_\pi) = |E(G_\pi)|$ and hence $\ell(\pi) = \ell'(\pi)$. Conversely, if $\ell(\pi) = \ell'(\pi)$, then the first inequality of Theorem~\ref{thm:permutation-chain}, due to Tenner and Vatter~\cite[Corollary 1.5]{TV}, forces $\beta(G_\pi) = 0$, so $G_\pi$ is a forest and $\pi$ avoids $321$ and $3412$.

\subsection*{Why $\graphrefl(G_\pi)$ can be smaller than $\ell'(\pi)$}

The inequality $\graphrefl(G_\pi) \le \ell'(\pi)$ of Theorem~\ref{thm:graph-length-at-most-reflection-length} can be strict, and the ways this happens can be grouped into two kinds. We take them in turn, and then give a family for which the difference $\ell'(\pi) - \graphrefl(G_\pi)$ grows without bound.

The first is that reflection length is not determined by the inversion graph at all. The \emph{diamond} $K_4 - e$ is the inversion graph of $3421$, $4231$, and $4312$, whose reflection lengths are $3$, $1$, and~$3$, respectively. Edge reflection length, by contrast, is an isomorphism invariant of the graph, so $\graphrefl(G_\sigma) = \graphrefl(G_\tau)$ whenever $G_\sigma \cong G_\tau$. These three permutations share an inversion graph, so any statistic that depends only on that graph must be constant on them; since $\ell'$ is not, no such statistic can equal $\ell'$, and in particular $\graphrefl$ does not. Section~\ref{sec:further-questions} takes up which permutation classes avoid this.

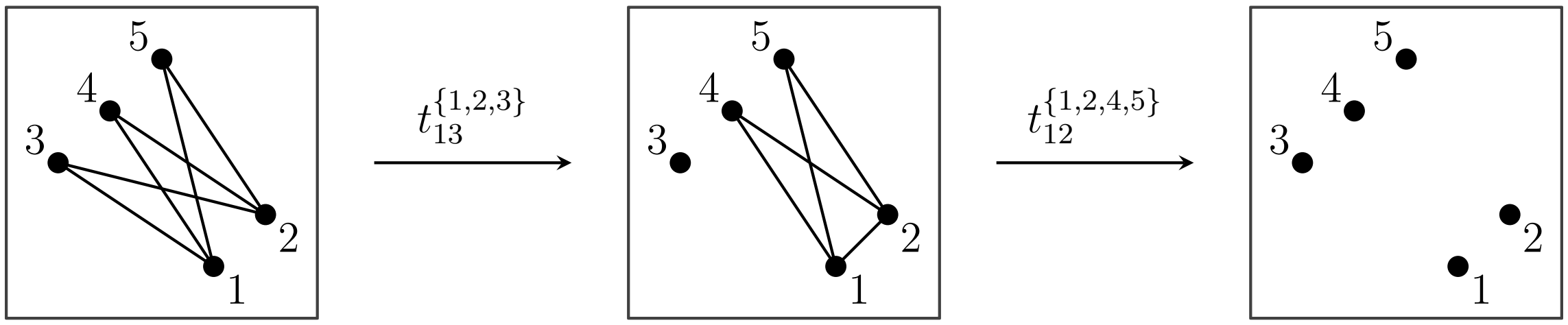
\begin{figure}[htbp]
\begin{center}
\begin{tikzpicture}[>=stealth]

\node (G1) at (-6, 0) {
	\begin{tikzpicture}[scale=0.5,baseline=(current bounding box.center)]
	\draw[darkgray,thick,rounded corners=0.01,line cap=round] (0,0) rectangle (6,6);
	\plotperm{3,4,5,1,2}
	\labelnode[135]{v3}{$3$}
	\labelnode[135]{v4}{$4$}
	\labelnode[135]{v5}{$5$}
	\labelnode[-45]{v1}{$1$}
	\labelnode[-45]{v2}{$2$}
    \drawjoin{v1,v2}{v3,v4,v5}
	\end{tikzpicture}
};

\node (G2) at (0, 0) {
	\begin{tikzpicture}[scale=0.5,baseline=(current bounding box.center)]
	\draw[darkgray,thick,rounded corners=0.01,line cap=round] (0,0) rectangle (6,6);
	\plotperm{3,4,5,1,2}
	\labelnode[135]{v3}{$3$}
	\labelnode[135]{v4}{$4$}
	\labelnode[135]{v5}{$5$}
	\labelnode[-45]{v1}{$1$}
	\labelnode[-45]{v2}{$2$}
    \drawjoin{v1,v2}{v2,v4,v5}
	\end{tikzpicture}
};

\node (G3) at (6, 0) {
	\begin{tikzpicture}[scale=0.5,baseline=(current bounding box.center)]
	\draw[darkgray,thick,rounded corners=0.01,line cap=round] (0,0) rectangle (6,6);
	\plotperm{3,4,5,1,2}
	\labelnode[135]{v3}{$3$}
	\labelnode[135]{v4}{$4$}
	\labelnode[135]{v5}{$5$}
	\labelnode[-45]{v1}{$1$}
	\labelnode[-45]{v2}{$2$}
	\end{tikzpicture}
};


\draw[->, thick, shorten <=14pt, shorten >=14pt]
  (G1) -- node[above, yshift=4pt] {$t_{13}^{\{1,2,3\}}$} (G2);

\draw[->, thick, shorten <=14pt, shorten >=14pt]
  (G2) -- node[above, yshift=4pt] {$t_{12}^{\{1,2,4,5\}}$} (G3);

\end{tikzpicture}
\end{center}
\caption{The graph $K_{3,2}$ (left) can be transformed to the edgeless graph $\overline{K}_5$ (right) with two edge reflections.}
\label{fig:K32-example}
\end{figure}

The second is that even when all the permutations realizing a graph share a reflection length, the edge reflection length can still be smaller, because of the freedom our definition of edge reflection allows in selecting the set~$X$. Take $K_{3,2}$, the inversion graph of precisely two permutations, $34512$ and $45123$, both of reflection length~$4$. Its edge reflection length is smaller: by Corollary~\ref{cor:exact-small-beta}, since $K_{3,2}$ is connected with $\beta(K_{3,2}) = 2$, we have $\graphrefl(K_{3,2}) = 6 - 4 = 2$. Figure~\ref{fig:K32-example} shows a pair of reflections achieving this, the first sending $K_{3,2}$ to a copy of $K_2 \vee \overline{K}_2$ and an isolated vertex. The freedom in choosing their sets~$X$ is what makes two reflections enough.

Finally, the difference $\ell'(\pi) - \graphrefl(G_\pi)$ can be arbitrarily large. The permutation $\pi = 34 \cdots n \, 21$ has inversion graph $K_2 \vee \overline{K}_{n-2}$, so $\graphrefl(G_\pi) = 1$ by Observation~\ref{obs:nested-triangles}. However, its reflection length grows with~$n$: when $n$ is even the permutation is a single $n$-cycle and $\ell'(\pi) = n - 1$, and when $n$ is odd it has two cycles and $\ell'(\pi) = n - 2$. Since $G_\pi$ is connected with $2n - 3$ edges, we have $\beta(G_\pi) = n - 2$, so for even~$n$ these permutations attain the upper bound of Corollary~\ref{cor:reflection-gap}. The same graph is the inversion graph of the reflection $t_{1n}$, whose reflection length is~$1$, so this is another example like the permutations of Figure~\ref{fig:diamond-realizers}, where reflection length is not determined by the inversion graph.

\section{Further questions}
\label{sec:further-questions}

We close with several questions arising from this work. The first and last concern edge reflections themselves: whether steps that increase the number of edges can help, and what the right analogue of edge reflections is for ordered graphs. The middle three concern the relation between inversion graphs and the permutations that realize them.

\subsection*{When does backtracking help?}

We only allow admissible edge reflections, so each edge reflection strictly decreases the number of edges. In this sense, a sequence of reflections that empties a graph is a sequence of ``forward'' steps, and we can ask what happens if we allow ``backward'' steps as well. Say that a step of a sequence \emph{backtracks} if it undoes an admissible reflection, passing from a graph~$H$ to a graph~$G$ for which $H = t_{uv}^X \circ G$ with $t_{uv}^X$ admissible at~$G$. Such a step increases the number of edges; is that ever helpful?

For permutations the answer is that it is not. Section~\ref{sec:sorting} produces, for every~$\pi$, a sequence of $\ell'(\pi)$ reflections sorting~$\pi$ in which every reflection is applied to an inversion, and hence decreases the Coxeter length. Dyer~\cite[Theorem 1.2]{dyer:on-minimal-leng:} proved the corresponding statement for an arbitrary Coxeter group: the reflection length of any element is achieved by a sequence of reflections each of which decreases the Coxeter length.

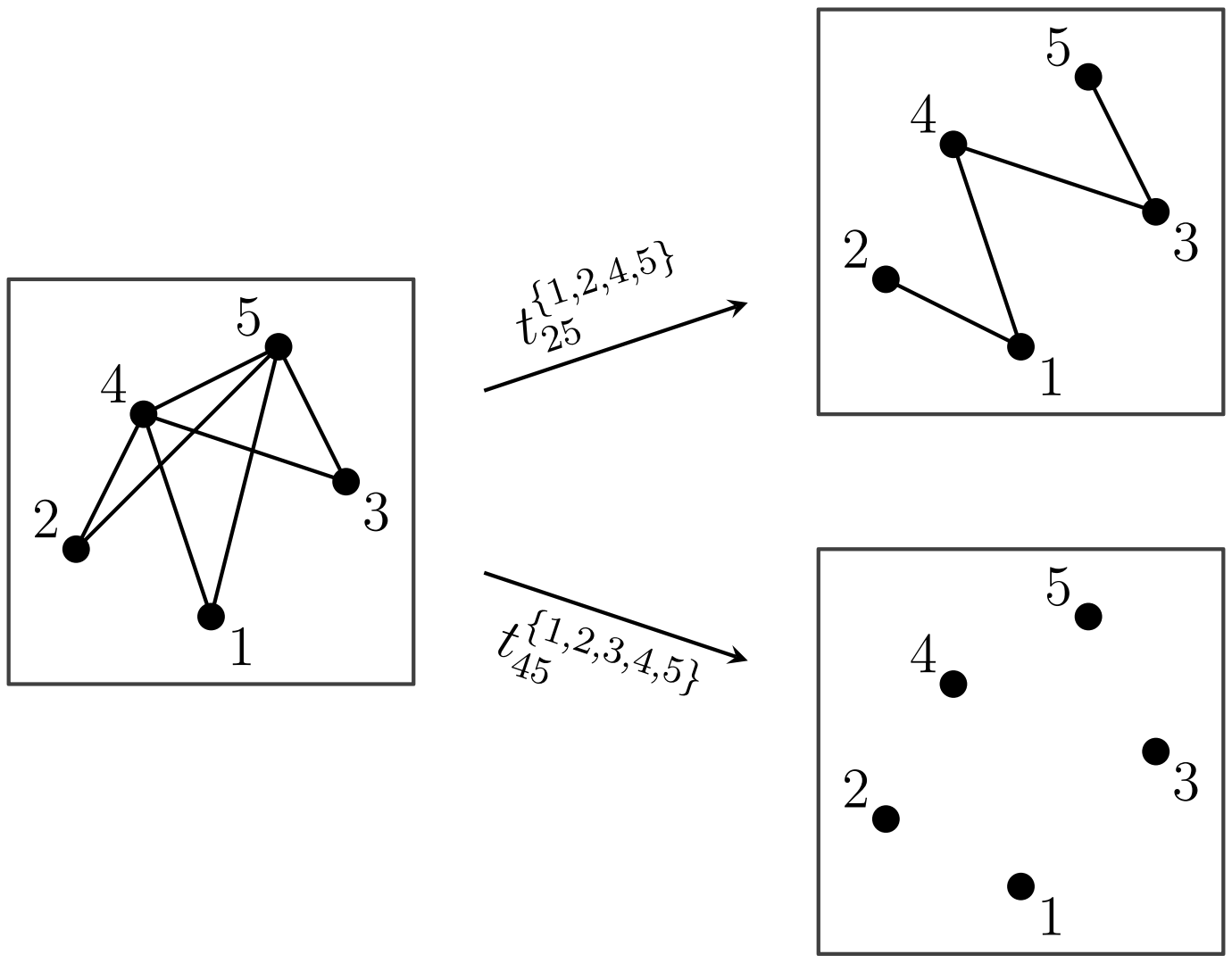
\begin{figure}[htbp]
\begin{center}
\begin{tikzpicture}[>=stealth]

\node (G2) at (0, 0) {
	\begin{tikzpicture}[scale=0.5,baseline=(current bounding box.center)]
	\draw[darkgray,thick,rounded corners=0.01,line cap=round] (0,0) rectangle (6,6);
	\plotperm{2,4,1,5,3}
	\labelnode[135]{v2}{$2$}
	\labelnode[135]{v4}{$4$}
	\labelnode[-45]{v1}{$1$}
	\labelnode[135]{v5}{$5$}
	\labelnode[-45]{v3}{$3$}
	\drawjoin{v4,v5}{v1,v3}
	\draw (v2)--(v4)--(v5)--(v2);
	\end{tikzpicture}
};

\node (G1) at (6, 2) {
	\begin{tikzpicture}[scale=0.5,baseline=(current bounding box.center)]
	\draw[darkgray,thick,rounded corners=0.01,line cap=round] (0,0) rectangle (6,6);
	\plotperm{2,4,1,5,3}
	\labelnode[135]{v2}{$2$}
	\labelnode[135]{v4}{$4$}
	\labelnode[-45]{v1}{$1$}
	\labelnode[135]{v5}{$5$}
	\labelnode[-45]{v3}{$3$}
	\draw (v2)--(v1)--(v4)--(v3)--(v5);
	\end{tikzpicture}
};

\node (G3) at (6, -2) {
	\begin{tikzpicture}[scale=0.5,baseline=(current bounding box.center)]
	\draw[darkgray,thick,rounded corners=0.01,line cap=round] (0,0) rectangle (6,6);
	\plotperm{2,4,1,5,3}
	\labelnode[135]{v2}{$2$}
	\labelnode[135]{v4}{$4$}
	\labelnode[-45]{v1}{$1$}
	\labelnode[135]{v5}{$5$}
	\labelnode[-45]{v3}{$3$}
	\end{tikzpicture}
};


\draw[->, thick, shorten <=14pt, shorten >=14pt]
  (G2) -- node[midway, above=2pt, sloped] {$t_{25}^{\{1,2,4,5\}}$} (G1);

\draw[->, thick, shorten <=14pt, shorten >=14pt]
  (G2) -- node[midway, below=2pt, sloped] {$t_{45}^{\{1,2,3,4,5\}}$} (G3);

\end{tikzpicture}
\end{center}
\caption{The graph $K_2 \vee \overline{K}_3$ (left) and two admissible edge reflections of it, one producing $G_{24153} \cong P_5$ (top right) and the other emptying it (bottom right). Emptying $P_5$ takes four edge reflections, but only two steps if we allow backtracking through $K_2 \vee \overline{K}_3$.}
\label{fig:P5-backtracking}
\end{figure}

However, for our notion of reflections in graphs, there is no analogue of Dyer's theorem; that is, backtracking can help. Consider, for example, the path $P_5$. This graph is a tree with four edges, so $\graphrefl(P_5) = 4$ by Corollary~\ref{cor:forest-edge-bound}. But, as shown in Figure~\ref{fig:P5-backtracking}, a single backtracking step to $K_2 \vee \overline{K}_3$, followed by one admissible reflection, empties the graph in two steps.

For permutations, the corresponding step saves nothing. Here $K_2 \vee \overline{K}_3$ is the inversion graph of $54123$, and $t_{25}$ sends $54123$ to $24153$, whose inversion graph is $P_5$, so the backward step is a legitimate permutation reflection. But $54123 = (1,5,3)(2,4)$ has two cycles, so $\ell'(54123) = 3$, and sorting $24153$ by way of $54123$ takes $1 + 3 = 4$ reflections, which is $\ell'(24153) = 4$, since $24153 = (1,2,4,5,3)$ is a single cycle. The advantage of backtracking comes from $\graphrefl(G_{54123}) = 1$, far below $\ell'(54123) = 3$.

\begin{question}
\label{q:backtracking}
For which graphs does backtracking provide an advantage when emptying the graph by edge reflections?
\end{question}

\subsection*{When is reflection length a graph invariant?}

As noted in Section~\ref{sec:permutations}, the diamond $K_4 - e$ is the inversion graph of $3421$, $4231$, and $4312$, of reflection lengths $3$, $1$, and~$3$, so reflection length is not determined by the inversion graph. It is, however, determined by the inversion graph within some permutation classes (sets of permutations closed downward in the pattern containment order; see Vatter~\cite{vatter:permutation-cla:} for a survey). For example, by Corollary~\ref{cor:forest-all-equal}, every $\pi \in \mathrm{Av}(321, 3412)$ satisfies $\ell'(\pi) = |E(G_\pi)|$. The question is which classes have this property.

\begin{question}
\label{q:graph-invariant}
For which permutation classes $\mathcal{C}$ is reflection length a graph invariant, meaning that $G_\sigma \cong G_\tau$ implies $\ell'(\sigma) = \ell'(\tau)$ for all $\sigma, \tau \in \mathcal{C}$?
\end{question}

Since $\ell'(\pi) = n - c(\pi)$, this asks for which classes the number of permutation cycles is determined by the inversion graph. The diamond already restricts the candidates. The three permutations with this inversion graph split by their number of cycles: $3421$ and $4312$ have one cycle, while $4231$ has three, so a class with this property cannot contain permutations from both groups. Because permutation classes are closed downward, every such class must avoid either both~$3421$ and~$4312$, or the single pattern~$4231$.

After searching all members through length~$10$, we have found no two permutations in $\mathrm{Av}(3421, 4312)$ that have isomorphic inversion graphs but different reflection lengths. At length~$10$, the class has $105{,}632$ members, realizing $12{,}128$ isomorphism classes of inversion graphs. This leads us to the following conjecture.

\begin{conjecture}
\label{conj:reflection-length-invariant}
The number of permutation cycles is a graph invariant for $\mathrm{Av}(3421, 4312)$.
\end{conjecture}

On the other hand, the class $\mathrm{Av}(4231)$ does not have the property. The four permutations
\[
	465132 \quad 465213 \quad 546132 \quad 546213
\]
all avoid $4231$ and share an inversion graph, but the first and last have cycle type $(2,2,2)$ while the other two are $6$-cycles, so their cycle counts differ (the corresponding reflection lengths are~$3$,~$5$,~$5$, and~$3$). Any subclass of $\mathrm{Av}(4231)$ on which the number of cycles is a graph invariant must therefore exclude both $465132$ and $546213$, or else both $465213$ and $546132$.

\subsection*{When is cycle type a graph invariant?}

A stronger property than cycle count invariance would be invariance of the full \emph{cycle type}, the multiset of cycle lengths of the permutation.

\begin{question}
\label{q:cycle-type}
For which permutation classes does $G_\sigma \cong G_\tau$ imply that $\sigma$ and $\tau$ have the same cycle type?
\end{question}

The \emph{layered permutations}, $\mathrm{Av}(231, 312)$, have this property, as we briefly explain. A layered permutation is a direct sum of decreasing permutations. If these have sizes $a_1, \ldots, a_k$, then each is a clique in the inversion graph, so the graph is the disjoint union $K_{a_1} \cup \cdots \cup K_{a_k}$. A decreasing permutation of size~$a$ consists of $\lfloor a/2 \rfloor$ disjoint transpositions, together with a fixed point when~$a$ is odd. The whole cycle type is thus determined by the multiset $\{a_1, \ldots, a_k\}$, so cycle type is a graph invariant on the layered permutations. We conjecture that this property extends to the class from Conjecture~\ref{conj:reflection-length-invariant}, which contains the class of layered permutations. As before, our evidence is the computer search through length~$10$.

\begin{conjecture}
\label{conj:cycle-type-invariant}
Cycle type is a graph invariant for $\mathrm{Av}(3421, 4312)$.
\end{conjecture}

Cycle type invariance is stronger than cycle count invariance, so it too fails in $\mathrm{Av}(4231)$, and it can fail even where cycle counts agree. For example, six members of $\mathrm{Av}(4231)$ have inversion graph $K_2 \vee (K_2 \cup \overline{K}_2)$: $346521$, $435621$, $651243$, and $652134$ have cycle type $(3,3)$, while $354621$ and $651324$ have cycle type $(4,2)$. All six have two cycles, so a subclass with cycle type invariance must omit some of these permutations.

\subsection*{\boldmath When does $\ell'(\pi) = \graphrefl(G_\pi)$?}

Corollary~\ref{cor:unicyclic-equality} shows that $\ell'(\pi) = \graphrefl(G_\pi)$ whenever every connected component of~$G_\pi$ has cyclomatic number at most one, and Corollary~\ref{cor:reflection-gap} bounds the difference by $\beta(G_\pi)$ in general. This difference can be arbitrarily large, as the permutations $34 \cdots n\, 21$ of Section~\ref{sec:permutations} show.

Returning to the layered permutations, both quantities follow from the sizes $a_1, \ldots, a_k$ of the decreasing summands. A decreasing permutation of size~$a$ has $\lceil a/2 \rceil$ cycles, so it contributes $a - \lceil a/2 \rceil = \lfloor a/2 \rfloor$ to the reflection length, and its inversion graph $K_a$ has $\graphrefl(K_a) = \lfloor a/2 \rfloor$ by Section~\ref{sec:lower-bounds}. Summing over the summands shows that the two totals agree.

\begin{observation}
\label{obs:layered-equality}
Every layered permutation satisfies $\ell'(\pi) = \graphrefl(G_\pi)$.
\end{observation}

A characterization by classical pattern avoidance is not possible, because the permutations satisfying $\ell'(\pi) = \graphrefl(G_\pi)$ do not form a permutation class. Equality holds for $35241 = (1,3,2,5)(4)$: it has two cycles, so $\ell'(35241) = 3$, and its inversion graph also has edge reflection length $3$ (the graph is connected with seven edges on five vertices, so $\beta(G_{35241}) = 3$, and Theorem~\ref{thm:sandwich-improved} with parity leaves edge reflection length $7 - 6 = 1$ or $7 - 4 = 3$; $1$ is impossible by Observation~\ref{obs:nested-triangles} because the graph is not of the form $K_2 \vee \overline{K}_m$). Yet $35241$ contains $3421$, for which $\ell'(3421) = 3$ while $\graphrefl(G_{3421}) = 1$.

This question, like the two before it, asks how much of a permutation's cycle structure is visible in its inversion graph, which records only the inversions. Classical pattern avoidance cannot answer it, as we saw, but a richer notion of pattern might. Berman and Tenner~\cite{berman-tenner} introduced ``arrow patterns,'' which augment a classical pattern with information from the fundamental bijection, and used them to resolve a problem of Petersen and Tenner~\cite{petersen:the-depth-of-a-:} on shallow permutations; the arrow patterns they used can also be expressed with mesh patterns~\cite{branden:mesh-patterns-a:,tenner:mesh-patterns-w:}. Either notion may be the right tool here.

\begin{question}
\label{q:equality}
For which permutations does $\ell'(\pi) = \graphrefl(G_\pi)$?
\end{question}

\subsection*{What about ordered graphs?}

We arrived at edge reflections by discarding the notion of ``horizontally between'' and letting the vertices playing that role be chosen freely. An ordered graph, whose vertices are linearly ordered, supplies such a notion of its own: for an edge $uv$ with $u < v$, the vertices $w$ with $u < w < v$ lie between its endpoints. This suggests a notion of reflection closer to permutation reflections, which we leave as a question.

\begin{question}
\label{q:ordered-graphs}
What is the right analogue of edge reflections for ordered graphs, and what can be said about it?
\end{question}

\bibliographystyle{abbrv}
\bibliography{references}

\end{document}